\documentclass[aap,citesort,MSNbibl,dvips]{arximspdf}
\usepackage{mathbh}

\doi{10.1214/09-AAP647}
\volume{20}
\issue{6}
\pubyear{2010}
\firstpage{1989}
\lastpage{2021}

\makeatletter

\newcommand{\eqref}[1]{(\ref{#1})}

\newcommand{\eps}{\varepsilon}

\newcommand{\ssup}[1]{{({#1})}}

\newcommand{\R}{\mathbb{R}}

\newcommand{\N}{\mathbb{N}}

\newcommand{\prob}{\mathbb{P}}

\newcommand{\me}{\mathbb{E}}

\renewcommand{\P}{\mathbb{P}}

\newcommand{\one}{\mathbh{1}}

\newcommand{\tnu}{\tilde{\nu}}

\newcommand{\emm}{{\mathfrak m}}
\newcommand{\cl}{\operatorname{cl}}
\newcommand{\interior}{\operatorname{int}}

\newcommand{\skria}{{\mathcal A}}
\newcommand{\skrib}{{\mathcal B}}
\newcommand{\skric}{{\mathcal C}}

\newcommand{\skrie}{{\mathcal E}}

\newcommand{\skrig}{{\mathcal G}}
\newcommand{\skrih}{{\mathcal H}}

\newcommand{\skrik}{{\mathcal K}}

\newcommand{\skrim}{{\mathcal M}}

\newcommand{\skrin}{{\mathcal N}}
\newcommand{\skrip}{{\mathcal P}}

\newcommand{\skris}{{\mathcal S}}

\newcommand{\skrix}{{\mathcal X}}
\newcommand{\skriy}{{\mathcal Y}}

\newcommand{\sfrac}[2]{\frac{#1}{#2}}

\newtheorem{theorem}{Theorem}[section]
\newtheorem{lemma}[theorem]{Lemma}
\newtheorem{cor}[theorem]{Corollary}

\newproclaim{remark}{Remark}
\newproclaim{example}{Example}

\makeatother

\begin{document}
\begin{frontmatter}

\title{Large deviation principles for empirical measures of
colored random graphs}
\runtitle{Large deviation principles for colored random graphs}

\begin{aug}
\author[A]{\fnms{Kwabena} \snm{Doku-Amponsah}\ead[label=e1]{kdoku@ug.edu.gh}} and
\author[B]{\fnms{Peter} \snm{M\"{o}rters}\corref{}\thanksref{t1}\ead[label=e2]{maspm@bath.ac.uk}}
\runauthor{K. Doku-Amponsah and P. M\"{o}rters}
\affiliation{University of Ghana and University of Bath}
\address[A]{Department of Statistics\\
University of Ghana\\
P.O. Box LG 115\\
Legon-Accra\\
Ghana\\
\printead{e1}} 
\address[B]{Department of Mathematical Sciences\\
University of Bath\\
Claverton Down\\
Bath BA2 7AY\\
United Kingdom\\
\printead{e2}}
\end{aug}

\thankstext{t1}{Supported by an Advanced Research
Fellowship of the EPSRC.}

\pdfauthor{Kwabena Doku-Amponsah, Peter Morters}

\received{\smonth{10} \syear{2006}}
\revised{\smonth{7} \syear{2009}}

%
\begin{abstract}
For any finite colored graph
we define the \textit{empirical neighborhood measure}, which counts the
number of vertices of a given color connected to a given number of vertices
of each color, and the \textit{empirical pair measure}, which
counts the number of edges connecting each pair of colors.
For a class of models of sparse colored random graphs, we
prove large deviation principles for these empirical measures in the
weak topology. The rate functions governing our large deviation principles
can be expressed explicitly in terms of relative entropies.
We derive a large deviation principle for the degree distribution of
Erd\H{o}s--R\'enyi graphs near criticality.
\end{abstract}

%

\end{frontmatter}

\section{Introduction}

In this paper we study a random graph model where each vertex of the
graph carries a
random symbol, spin or color. The easiest model of this kind is that of
an Erd\H os--R\'enyi
graph where additionally each vertex is equipped with an independently
chosen color.
The more general models of colored random graphs we consider here allow
for a dependence between color
and connectivity of the vertices.

With each colored graph we associate its \textit{empirical neighborhood
measure}, which
records the number of vertices of a given color with a given number of
adjacent vertices
of each color. From this quantity one can derive a host of important
characteristics of the colored
graph, like its degree distribution, the number of edges linking two
given colors,
or the number of isolated vertices of any color. The aim of this paper
is to derive
a large deviation principle for the empirical neighborhood measure.

To be more specific about our model, we consider colored random graphs
constructed as follows: in the first step each of $n$ fixed vertices
independently
gets a color, chosen according to some law $\mu$ on the finite set
$\skrix$ of colors.
In the second step we connect any pair of vertices independently with a
probability
$p(a,b)$ depending on the colors $a,b\in\skrix$ of the two vertices.
This model, which comprises the simple Erd\H os--R\'enyi graph with independent
colors as a special case, was introduced by Penman in his thesis \cite{Pe98}
(see \cite{CP03} for an exposition) and rediscovered later by
S\"oderberg \cite{So02}.
It is also a special case of the inhomogeneous random graphs studied
in \cite{BJR07,vdH09}.

Our main concern in this paper are asymptotic results when the
number $n$ of vertices goes to infinity, while the connection
probabilities are of order $1/n$. This leads to an average number of
edges of order $n$, the \textit{near critical} or \textit{sparse} case.
Our methods also allow the study of the sub- and supercritical
regimes. Some results on these cases are discussed in \cite{DA06}.

Recall that a \textit{rate function} is a nonconstant, lower
semicontinuous function
$I$ from a polish space $\skrim$ into $[0,\infty]$, it is called \textit{good}
if the level sets $\{ I(m)\le x\}$ are compact for every $x\in[0,\infty)$.
A functional $M$ from the set of finite colored graphs to $\skrim$ is
said to
satisfy a \textit{large deviation principle} with rate function $I$ if,
for all
Borel sets $B\subset\skrim$,
\begin{eqnarray*}
-\inf_{m\in\interior B} I(m) & \le& \liminf_{n \to\infty} \sfrac1n
\log\prob_n \{ M(X) \in B \} \\
& \le& \limsup_{n \to\infty} \sfrac1n \log\prob_n \{ M(X) \in B \}\\
& \le& -\inf_{m\in\cl B} I(m),
\end{eqnarray*}
where $X$ under $\prob_n$ is a colored random graph with $n$ vertices
and int $B$ and $\operatorname{cl}B$ refer to
the interior, respectively, closure, of the set $B$.

Apart from the empirical neighborhood measure defined above, we also consider
the \textit{empirical pair measure}, which counts the number of edges
connecting any given pair of colors,
and the \textit{empirical color measure}, which simply counts the number
of vertices of any given color.
The main result of this paper is a joint large deviation principle
for the empirical neighborhood measure and the empirical pair measure
of a colored random graph
in the weak topology; see Theorem \ref{main2}. In the course of the
proof of this principle, two
further interesting large deviation principles are established: a large
deviation principle for the empirical
neighborhood measure conditioned to have a given empirical pair and
color measure
(see Theorem \ref{main3}) and a joint large deviation principle for the
empirical color measure and the
empirical pair measure; see Theorem \ref{main}(b). For all these
principles we obtain a
completely explicit rate function given in terms of relative
entropies.

Our motivation for this project is twofold. On the one hand, one may consider
the colored random graphs as a very simple random model of
\textit{networked data}.
The data is described as a text of fixed length, consisting of words
chosen from
a finite dictionary, together with a random number of unoriented edges
or links connecting
the words. Large deviation results for the empirical neighborhood
measure permit
the calculation of the asymptotic number of bits needed to transmit a
large amount
of such data with arbitrarily small error probability; see \cite{DA06}
where this idea is followed up.

On the other hand, we are working toward understanding simple models of
\textit{statistical mechanics} defined on random graphs. Here, typically,
the colors of the vertices are
interpreted as spins, taken from a finite set of possibilities, and the
Hamiltonian of
the system is an integral of some function with respect to the
empirical neighborhood
measure. As a very simple example we provide the annealed asymptotics of
the random partition function for the Ising model on an
Erd\H{o}s--R\'enyi graph,
as the graph size goes to infinity.

As a more substantial example, we consider the Erd\H{o}s--R\'enyi graph
model on $n$ vertices,
where edges are inserted with probability $p_n\in[0,1]$ independently
for any pair of vertices. We assume that $np_n \to c\in(0,\infty)$.
From our main result we derive a large
deviation principle for the degree distribution; see Corollary \ref{ERdd}.
This example seems to be new in this explicit form.

\section{Statement of the results}\label{sec2}

Let $V$ be a fixed set of $n$ vertices, say $V=\{1,\ldots,n\}$
and denote by $\skrig_n$ the set of all (simple) graphs with
vertex set $V=\{1,\ldots,n\}$ and edge set
$E\subset\skrie:= \{(u,v)\in V\times V \dvtx u<v \}$, where
the \textit{formal} ordering of edges is introduced as a means to
describe simply \textit{unordered} edges. Note that for all $n$, we
have $0\le|E|\le\sfrac{1}{2} n(n-1)$. Let $\skrix$ be a
finite alphabet or color set $\skrix$ and denote by
$\skrig_n(\skrix)$ be the set of all colored graphs with color set
$\skrix$ and $n$ vertices.

Given a symmetric function $p_n\dvtx\skrix\times\skrix\rightarrow
[0,1]$ and a probability measure $\mu$ on $\skrix$ we may define the
\textit{randomly colored random graph} or simply \textit{colored random
graph} $X$ with $n$ vertices as follows: assign to each vertex $v\in
V$ color $X(v)$ independently according to the \textit{color law}
$\mu$. Given the colors, we connect any two vertices $u,v\in V$,
independently of everything else, with \textit{connection
probability} $p_n(X(u),X(v))$.
We always consider $X=((X(v) \dvtx v\in V),E)$ under the joint law of
graph and
color and interpret $X$ as colored random graph.

We are interested in the properties of the randomly colored graphs for
large $n$ in the \textit{sparse} or \textit{near critical} case, that
is, we
assume that the connection probabilities satisfy $n p_n(a,b) \to C(a,b)$
for all $a,b\in\skrix$, where $C\dvtx\skrix\times\skrix\rightarrow
[0,\infty)$ is a symmetric function, which is not identically equal to zero.

To fix some notation, for any finite or countable set $\skriy$ we
denote by
$\skrim(\skriy)$ the space of probability measures, and by
$\tilde\skrim(\skriy)$ the space of finite measures on $\skriy$, both
endowed with the weak topology. When applying $\nu\in\tilde\skrim(\skriy)$
to some function $g\dvtx\skriy\to\R$ we use the scalar-product notation
\[
\langle\nu, g\rangle:= \sum_{y\in\skriy} \nu(y) g(y)
\]
and denote by $\|\nu\|$ its total mass. Further, if $\mu\in\tilde\skrim
(\skriy)$ and
$\nu\ll\mu$ we denote by
\[
H(\nu\| \mu)=\sum_{y\in\skriy} \nu(y) \log\biggl(\sfrac{\nu(y)}{\mu(y)} \biggr)
\]
the \textit{relative entropy} of $\nu$ with respect to $\mu$. We set
$H(\nu\| \mu)=\infty$ if $\nu\not\ll\mu$. By $\skrin(\skriy)$ we
denote the space of
counting measures on $\skriy$, that is, those measures taking\vspace*{1pt} values in
$\N\cup\{0\}$,
endowed with the discrete topology. Finally, we denote by $\tilde\skrim
_*(\skriy\times\skriy)$
the subspace of symmetric measures in $\tilde\skrim(\skriy\times\skriy
)$.%

With any colored graph $X=((X(v) \dvtx v\in V),E)$ with $n$ vertices
we associate a probability measure, the \textit{empirical color
measure} $L^1\in\skrim(\skrix)$, by
\[
L^{1}(a):=\frac{1}{n}\sum_{v\in V}\delta_{X(v)}(a) \qquad\mbox{for $a\in
\skrix$,}
\]
and a symmetric finite measure, the \textit{empirical pair measure}
$L^{2}\in\tilde\skrim_*(\skrix\times\skrix)$, by
\[
L^{2}(a,b):=\frac{1}{n}\sum_{(u,v)\in E}\bigl[\delta_{(X(v), X(u))}+
\delta_{(X(u), X(v))}\bigr](a,b) \qquad\mbox{for $a,b\in\skrix$.}
\]
The total mass $\|L^2\|$ of the empirical pair measure is
$2|E|/n$. Finally we define a further probability measure, the
\textit{empirical neighborhood measure}
$M\in\skrim(\skrix\times\skrin(\skrix))$, by
\[
M(a,\ell):=\frac{1}{n}\sum_{v\in V}\delta_{(X(v),L(v))}(a,\ell)\qquad
\mbox{for $(a,\ell)\in\skrix\times\skrin(\skrix)$,}
\]
where
$L(v)=(l^{v}(b), b\in\skrix)$ and $l^{v}(b)$ is the number of
vertices of color $b$ connected to vertex $v$. For every
$\nu\in\skrim(\skrix\times\skrin(\skrix))$ let $\nu_1$ and $\nu_2$
be the $\skrix$-marginal and the $\skrin(\skrix)$-marginal of the
measure $\nu$, respectively. Moreover, we define a measure
$\langle\nu(\cdot,\ell), \ell(\cdot)\rangle\in\tilde\skrim(\skrix\times
\skrix)$
by
\[
\langle\nu(\cdot,\ell), \ell(\cdot)\rangle(a,b):=
\sum_{\ell\in\skrin(\skrix)}\nu(a,\ell)\ell(b) \qquad\mbox{for $a,b\in\skrix$.}
\]
Define the function $\Phi\dvtx\skrim(\skrix\times\skrin(\skrix)) \to
\skrim(\skrix) \times\tilde\skrim(\skrix\times\skrix)$ by
$\Phi(\nu)=(\nu_1$,\break $\langle\nu(\cdot,\ell), \ell(\cdot)\rangle)$,
and observe that $\Phi(M)=(L^1, L^2)$, if these quantities are defined
as empirical
neighborhood, color, and pair measures of a colored graph. Note that
while the
\textit{first} component of $\Phi$ is a continuous function,
the \textit{second} component is \textit{discontinuous} in the weak topology.

To formulate the large deviation principle, we call a pair of measures
$(\varpi,\nu)\in\tilde{\skrim}(\skrix\times\skrix)\times\skrim(\skrix
\times\skrin(\skrix))$
\textit{sub-consistent} if
%
%
\begin{equation}\label{consistent}
\langle\nu(\cdot,\ell), \ell(\cdot)\rangle(a,b) \le\varpi(a,b)\qquad
\mbox{for all $a,b\in\skrix$,}
\end{equation}
and \textit{consistent} if equality holds in \eqref{consistent}.
Observe that, if $\nu$ is the empirical neighborhood measure and
$\varpi$ the empirical pair measure of a colored graph, $(\varpi, \nu)$ is
consistent and both sides in \eqref{consistent} represent
\[
\sfrac{1}{n} \bigl(1+\one_{\{a=b\}}\bigr) \sharp\{ \mbox{edges between
vertices of
colors $a$ and $b$} \}.
\]
For a measure $\varpi\in\tilde\skrim_*(\skrix\times\skrix)$ and a measure
$\omega\in\skrim(\skrix)$, define
\[
{{\mathfrak H}_C}(\varpi \| \omega):=
H (\varpi\| C\omega\otimes\omega)+\| C\omega\otimes\omega
\| -\|\varpi\| ,
\]
where the measure
$C\omega\otimes\omega\in\tilde\skrim(\skrix\times\skrix)$ is defined
by $C\omega\otimes\omega(a,b)=C(a,b)\*\omega(a)\omega(b)$ for
$a,b\in\skrix$. It is not hard to see that $\mathfrak
H_C(\varpi\| \omega)\ge0$ and equality holds if and only if
$\varpi=C\omega\otimes\omega$ (see Lemma \ref{randomg.Vrate}). For
every $(\varpi,\nu)\in\tilde\skrim_*(\skrix\times\skrix) \times
\skrim(\skrix\times\skrin(\skrix))$ define a probability measure
$Q=Q[\varpi,\nu]$ on $\skrix\times\skrin(\skrix)$ by
%
%
\begin{eqnarray}\label{Poissonlimit}
Q(a , \ell):=\nu_{1}(a)\prod_{b\in\skrix}
e^{-{\varpi(a,b)}/{\nu_1(a)}} \frac{1}{\ell(b)!}
\biggl(\frac{\varpi(a,b)}{\nu_1(a)} \biggr)^{\ell(b)} \nonumber\\[-8pt]\\[-8pt]
\eqntext{\mbox{for
$a\in\skrix$, $\ell\in\skrin(\skrix)$.}}
\end{eqnarray}
We have now set the stage to state our principal theorem, the large
deviation principle for the empirical pair measure and the empirical
neighborhood measure.
\begin{theorem}\label{main2}
Suppose that $X$ is a colored random graph with color law
$\mu$ and connection probabilities $ p_n\dvtx\skrix\times\skrix
\rightarrow[0,1]$
satisfying $n p_n(a,b) \to C(a,b)$ for some symmetric function
$C\dvtx\skrix\times\skrix\rightarrow[0,\infty)$ not identical to zero.
Then, as $n\rightarrow\infty$, the pair $(L^2, M)$ satisfies a large
deviation principle in
$\tilde{\skrim}_*(\skrix\times\skrix)\times\skrim(\skrix\times\skrin
(\skrix))$
with good rate function
\[
J(\varpi,\nu)=\cases{H(\nu\| Q)+H(\nu_1 \| \mu)+
\sfrac{1}{2} {{\mathfrak H}_C}(\varpi\| \nu_1 ), &\quad if $(\varpi,\nu)$
sub-consistent,\cr
\infty, &\quad otherwise.}
\]
\end{theorem}
\begin{remark}
The rate function can be interpreted as follows: $J(\varpi,\nu)$
represents the cost of obtaining an empirical pair
measure $\varpi$ and an empirical neighborhood measure $\nu$. This cost
is divided into three sub-costs:
\begin{longlist}
\item$H(\nu_1 \| \mu)$ represents the cost of obtaining the
empirical color measure $\nu_1$, this
cost is nonnegative and vanishes iff $\nu_1=\mu$.

\item$\sfrac{1}{2}{\mathfrak{H}_C}(\varpi\| \nu_1 )$
represent the cost of obtaining an empirical pair
measure $\varpi$ if the empirical color measure is $\nu_1$, again this
cost is nonnegative and vanishes iff
$\varpi= C \nu_1 \otimes\nu_1$.

\item$H(\nu\| Q)$ represents the cost of obtaining an
empirical neighborhood
measure $\nu$ if the empirical color measure is $\nu_1$ and the
empirical pair measure is
$\varpi$, this cost is nonnegative and vanishes iff $\nu=Q$.
\end{longlist}
Consequently, $J(\varpi,\nu)$ is nonnegative and vanishes if and only
if $\varpi= C \mu\otimes\mu$ and
\[
\nu(a,\ell)=\mu(a)\prod_{b\in\skrix}e^{-C(a,b)\mu(b)} \frac{(C(a,b)\mu
(b))^{\ell(b)}}{\ell(b)!} \qquad\mbox{for
all $(a,\ell)\in\skrix\times\skrin(\skrix)$}.
\]
This is the law of a pair $(a,\ell)$ where $a$ is distributed according
to $\mu$ and, given the value of $a$, the
random variables $\ell(b)$ are independently Poisson distributed with
parameter $C(a,b)\mu(b)$. Hence, as $n\to\infty$,
the random measure $M(X)$ converges in probability to this law.
\end{remark}
\begin{remark}
Our large deviation principle implies individual large deviation
principles for
the measures $L^2$ and $M$ by contraction; see \cite{DZ98},
Theorem 4.2.1. Note that,
by the discontinuity of $\Phi$, the functional relationship $L^2=\Phi
_2(M)$ may break down in the limit,
and hence the rate function may be finite on pairs which are not
consistent. We have not been able
to extend the large deviation principle to a \textit{stronger topology}
in which
$\Phi$ is continuous, as this leads to considerable compactness
problems; see \cite{ES02,DMS05}
for discussions of some of the problems and opportunities arising when
extending large deviation
principles to stronger topologies.
\end{remark}

As usual, the \textit{degree distribution} $D \in\skrim(\N\cup\{0\})$
of a graph with empirical neighborhood measure $M$
is defined by
\[
D(k)= \sum_{a\in\skrix}\sum_ {\ell\in\skrin(\skrix)} \delta_k
\biggl({ \sum_{b\in\skrix}} \ell(b) \biggr)
M(a,\ell) \qquad\mbox{for $k\in\N\cup\{0\}$,}
\]
that is, $D(k)$ is the proportion of vertices in the graph with degree $k$.
As the degree distribution $D$ is a continuous function of $M$,
Theorem \ref{main2}
and the contraction principle imply a large deviation principle for
$D$. For a classical
Erd\H{o}s--R\'enyi graph the rate function takes on a particularly
simple form (see
Section \ref{proofcorollaries} for details).
\begin{cor}\label{ERdd}
Suppose $D$ is the degree distribution of an Erd\H{o}s--R\'enyi graph
with connection probability
$p_n\in[0,1]$ satisfying $n p_n \to c \in(0,\infty)$. Then $D$
satisfies a large
deviation principle, as $n\to\infty$, in the space $\skrim(\N\cup\{0\}
)$ with good rate function
%
%
\begin{equation}\label{randomg.ratedeg}\qquad
\delta(d)= \cases{
\displaystyle\sfrac12 x \log\biggl( \sfrac xc \biggr) - \sfrac
12 x +\sfrac c2+H (d \| q_{x}), &\quad if $\langle d\rangle\le c$,\vspace*{2pt}\cr
\displaystyle\sfrac12 \langle d \rangle
\log\biggl( \sfrac{\langle d \rangle}{c} \biggr)
- \sfrac12 \langle d \rangle+ \sfrac c2 +H \bigl(d \| q_{\langle
d\rangle}\bigr),
&\quad if $c< \langle d\rangle< \infty$, \vspace*{2pt}\cr
\infty, &\quad if $\langle d\rangle= \infty$,}
\end{equation}
where, in the case $\langle d\rangle\le c$, the value $x=x(d)$ is the
unique solution of
\[
x=ce^{-2 (1-{\langle
d\rangle}/{x} )},
\]
and where $q_{\lambda}$ is a Poisson distribution with
parameter $\lambda$ and $\langle d\rangle:= \sum_{m=0}^{\infty}md(m)$.
\end{cor}
\begin{remark}
On probability measures $d$ with mean $c$ the rate simplifies to the
relative entropy
of $d$ with respect to the Poisson distribution of the same mean.
In~\cite{BGL02}, Theorem 7.1, a large deviation principle for the
degree distribution
is formulated for this situation, albeit with a rather implicitly
defined rate function.
Moreover, the proof given there contains a serious gap: the exponential
equivalence stated
in \cite{BGL02}, Lemma 7.2, is not proved there and we conjecture that
it does not hold.
\end{remark}
\begin{remark}
O'Connell \cite{OC98} provides further large deviation principles for sparse
Erd\H{o}s--R\'enyi graphs. His focus is on the size of the biggest
component, and
he also studies the number of isolated vertices. A large deviation
principle for the latter
is also a consequence of our corollary.
\end{remark}

%

We now state the two large deviation results, Theorems \ref{main}(b)
and \ref{main3}, which are
the main ingredients for our proof of Theorem \ref{main2}, but are also
of independent interest.
The first of these is a joint large deviation principle for the
empirical color measure $L^{1}$ and the
empirical pair measure $L^{2}$, the second a large deviation principle
for the empirical
neighborhood measure $M$ given $L^1$ and $L^2$.

For any $n\in\N$ we define
\begin{eqnarray*}
\skrim_n(\skrix) & := & \{ \omega\in\skrim(\skrix) \dvtx n\omega(a)
\in\N\mbox{ for all } a\in\skrix\},\\
\tilde\skrim_{*,n}(\skrix\times\skrix) & := & \biggl\{ \varpi\in\tilde\skrim
_*(\skrix\times\skrix)
\dvtx \sfrac n{1+\one\{a=b\}} \varpi(a,b) \in\N\\
&&\hspace*{140.41pt}\mbox{for all } a,b\in\skrix\biggr\}.
\end{eqnarray*}

\begin{theorem}\label{main}
Suppose that $X$ is a colored random graph with color law
$\mu$ and connection probabilities satisfying $n p_n(a,b) \to C(a,b)$
for some symmetric function $C\dvtx\skrix\times\skrix\rightarrow
[0,\infty)$
not identical to zero.

\textup{(a)} Suppose the sequence $\omega_n\in\skrim_n(\skrix)$ converges
to a
limit $\omega\in\skrim(\skrix)$. Then, as $n\rightarrow\infty$, conditional
on the event $\{L^1=\omega_n\}$ the empirical pair\vspace*{1pt}
measure $L^2$ of $X$ satisfies a large deviation principle on the
space $\tilde\skrim_*(\skrix\times\skrix)$ with good rate function
%
%
\begin{equation}\label{randomg.rateLDprobb}
I_{\omega}(\varpi)=\tfrac{1}{2}{\mathfrak{H}_C}(\varpi\| \omega) .
\end{equation}

\textup{(b)} As $n\rightarrow\infty$, the pair $(L^1,L^2)$ satisfies a
large deviation principle in
$\skrim(\skrix)\times\tilde{\skrim}_*(\skrix\times\skrix)$ with good
rate function
%
%
\begin{equation}\label{randomg.rateL2L1}
I(\omega,\varpi)=H(\omega\| \mu)+\tfrac{1}{2}{\mathfrak{H}_C}(\varpi
\| \omega) .
\end{equation}
\end{theorem}
%
%
\begin{example}\label{exa}
We look at the Erd\H os--R\'enyi graph with connection probabilities
$p_n$ satisfying
$np_n \to c\in(0,\infty)$ and study the random
partition function for the \textit{Ising model} on the graph, which is
defined as
\begin{eqnarray}
Z(\beta) := \sum_{{\eta\in\{-1,+1\}}^{V}} \exp\biggl( \beta\sum_{(u,v)\in
E} \eta(u)\eta(v) \biggr)\nonumber\\
\eqntext{\mbox{for the inverse temperature $\beta>0$.}}
\end{eqnarray}
Denoting by ${\sf E}$ expectation with respect to the graph, we note that
\[
{\sf E}Z(\beta)=2^n \me\exp\biggl( n \sfrac\beta2 \int xy
L^2(dx \, dy) \biggr) ,
\]
where $\me$ is expectation with respect to the graph randomly
colored using independent colors chosen uniformly from
$\skrix=\{-1,1\}$. Then Varadhan's lemma (see,
e.g., \cite{dH00}, III.3, Theorem \ref{main}(b)) gives
%
%
\begin{eqnarray}
\label{happy}
&&\lim_{n\to\infty}  \sfrac1n \log{\sf E} Z(\beta) \nonumber\\
&&\qquad = \log2 + \sup\biggl\{ \sfrac\beta2 \int xy \varpi(dx \,dy)-
I(\omega, \varpi) \dvtx
\omega\in\skrim(\skrix),\nonumber\\
&&\hspace*{148.1pt} \varpi\in\skrim_*(\{-1,1\} \times\{-1,1\})
\biggr\} \\
&&\qquad = \sup\biggl\{ \sfrac\beta2 \bigl( \varpi(\Delta) - \varpi(\Delta^{\mathrm{c}}) \bigr)
-x \log(x)\nonumber\\
&&\hspace*{52.9pt}{}-(1-x)\log(1-x) - \sfrac12 \bigl( H(\varpi\| \omega_x) + c -
\|\varpi\| \bigr) \biggr\} ,\nonumber
\end{eqnarray}
where $\Delta$ is the diagonal in $\{-1,1\} \times\{-1,1\}$, and the
supremum is over all $x\in[0,1]$
and $\varpi\in\skrim_*(\{-1,1\} \times\{-1,1\})$,
and the measure $\omega_x\in\tilde\skrim_*(\{-1,1\}\times\{-1,1\})$ is
defined by
\[
\omega_x(i,j)=cx^{(2+i+j)/2}(1-x)^{(2-i-j)/2} \qquad\mbox{for } i,j\in\{
-1,1\} .
\]
Note that the last expression in \eqref{happy} is an optimisation
problem in only four real variables.
\end{example}

We obtain from Theorem \ref{main}(b) the following corollary (see
Section \ref{proofcorollaries} for details).
\begin{cor}\label{randomg.L1E}
Suppose that $X$ is a colored random graph with color law $\mu$ and
connection probabilities satisfying
$n p_n(a,b) \to C(a,b)$ for some symmetric function
$C\dvtx\skrix\times\skrix\rightarrow[0,\infty)$ not identical to zero.
Then, as $n\rightarrow\infty$, the number
of edges per vertex $|E|/n$ satisfies a large deviation principle in
$[0,\infty)$ with good rate function
\[
\zeta(x)= x \log x - x+ \inf_{y>0} \biggl\{ \psi(y) - x \log\biggl(\sfrac12 y\biggr) +
\sfrac12 y \biggr\},
\]
where $\psi(y) = \inf H(\omega \| \mu)$ over all probability vectors
$\omega$ with $\omega^{T} C \omega=y$.
\end{cor}
\begin{remark}
In the Erd\H os--R\'enyi case $C(a,b)=c$ one obtains $\psi(y)=0$ for
$y=c$, and $\psi(y)=\infty$
otherwise. Hence\vspace*{1pt} $\zeta(x)= x \log x -x - x \log(\sfrac c2) + \sfrac
c2$, which is
the Cram\'er rate function for the Poisson distribution with
parameter $\frac c2$.
In \cite{BP03}, a large deviation principle for $|E|/n^2$ is proved for
colored random graphs with \textit{fixed} connection probabilities.
\end{remark}


For a given $\varpi\in\tilde\skrim_*(\skrix\times\skrix)$
and $\nu\in\skrim(\skrix\times\skrin(\skrix))$ we recall the
definition of
the measure $Q\in\skrim(\skrix\times\skrin(\skrix))$ from
\eqref{Poissonlimit}.
%
%
\begin{theorem}\label{randomg.LDprob}\label{main3}
Suppose $(\omega_n,\varpi_n)\in\skrim_n(\skrix)\times\tilde{\skrim
}_{*,n}(\skrix\times\skrix)$
converges to a limit $(\omega,\varpi)\in\skrim(\skrix)\times\tilde
{\skrim}_*(\skrix\times\skrix)$.
Let $X$ be a colored random graph with $n$ vertices conditioned on the
event $\{ \Phi(M)=(\omega_n,\varpi_n) \}$.
Then, as $n\rightarrow\infty$, the empirical neighborhood measure $M$
of $X$ satisfies a large deviation principle
in the space $\skrim(\skrix\times\skrin(\skrix))$ with good rate function
%
%
\begin{equation}\label{randomg.rateLDprob}
\tilde{J}_{(\omega,\varpi)}(\nu)=
\cases{H(\nu\| Q), &\quad if $(\varpi,{\nu})$ is sub-consistent and
${\nu}_1=\omega$,\cr
\infty, &\quad otherwise.}
\end{equation}
\end{theorem}



In the remainder of the paper we give the proofs of the results set
out so far. Section \ref{proofmain} is devoted to the proof of
Theorem \ref{main}(a), which uses the G\"{a}rtner--Ellis theorem.
By contrast, the proof of Theorem \ref{main3}, carried out in
Section \ref{proofmain3}, is based on nontrivial combinatorial
arguments combined with a partially randomized approximation
procedure. This approximation is the most demanding argument of the
paper and requires a fairly sophisticated technique. In
Section \ref{proofmain2} we first combine Sanov's theorem
\cite{DZ98}, Theorem 2.1.10 and Theorem \ref{main}(a) to obtain
Theorem~\ref{main}(b), and then Theorem \ref{main}(b) and
Theorem \ref{main3} to get Theorem \ref{main2}, using the setup and
result of Biggins \cite{Bi04} to ``mix'' the large deviation
principles. The paper concludes with the proofs of
Corollaries \ref{ERdd} and \ref{randomg.L1E}, which are given in
Section~\ref{proofcorollaries}. 

\section{\texorpdfstring{Proof of Theorem \protect\ref{main}\textup{(a)} by the
G\"{a}rtner--Ellis theorem}{Proof of Theorem 2.3(a) by the
G\"{a}rtner--Ellis theorem}}\label{proofmain}

Throughout this section we assume that
the sequence $\omega_n\in\skrim_n(\skrix)$ converges to
$\omega\in\skrim(\skrix)$. Let $\P\{ \cdot| L^1=\omega_n\}$ be the
law of colored random
graph $X$ with connection probabilities satisfying $n p_n(a,b) \to
C(a,b)$ conditioned on
the event $\{L^1=\omega_n\}$.
%
In the next lemma we verify the assumption of the G\"{a}rtner--Ellis theorem
\cite{DZ98}, Theorem 2.3.6. We denote by $\skric_2$ the space of symmetric
functions on $\skrix\times\skrix$.
\begin{lemma}\label{random.logg}
For each $g\in\skric_2$, we have
\[
\lim_{n\to\infty}\sfrac{1}{n}\log\me\bigl[e^{ n\langle
g, L^2\rangle}|L^1=\omega_{n} \bigr]=- \sfrac12 \langle
C\omega\otimes\omega, (1-e^{g})\rangle.
\]
\end{lemma}
\begin{pf}
Let $g\in\skric_2$. Observe that given the colors $a,b\in\skrix$ the
random variables
$nL^2(a,b)$ are binomial with parameters $n^2\omega_n(a)\omega
_n(b)-n\omega_n(a)\one_{\{a=b\}}$
and $p_n(a,b)$, and the variables $nL^2(a,b)$, for $\{a,b\}\subset\skrix
$ are independent.
Hence, we have that
\begin{eqnarray*}
&&\me\bigl[e^{ n\langle g, L^2\rangle} | L^1=\omega_{n} \bigr]
\\
&&\qquad=\prod_{\{a,b\}} \bigl(1-p_n(a,b)+p_n(a,b)e^{g(a,b)} \bigr)^{n^2\omega_n(a)\omega
_n(b)-n\omega_n(a)\one_{\{a=b\}} }.
\end{eqnarray*}
Now, for any $\eps>0$ and for large $n$ we have
\begin{eqnarray*}
\biggl(1-\sfrac{C(a,b)(1-e^{g(a,b)})+\eps}{n} \biggr)^n&\le&
\bigl(1-p_n(a,b)+p_n(a,b)e^{g(a,b)} \bigr)^{n}\\
&\le&\biggl(1-\sfrac{C(a,b)(1-e^{g(a,b)})-\eps}{n} \biggr)^n.
\end{eqnarray*}
Using Euler's formula and taking the product over all $\{a,b\}\subset
\skrix$ we obtain
\begin{eqnarray*}
-\sfrac12 \langle C\omega\otimes\omega, (1-e^{g})\rangle-\eps&\le&
\lim_{n\to\infty}\sfrac{1}{n}\log\me\bigl[e^{ n\langle
g, L^2\rangle}|L^1=\omega_{n} \bigr]\\
&\le&- \sfrac12 \langle
C\omega\otimes\omega, (1-e^{g})\rangle+\eps,
\end{eqnarray*}
and the result follows as $\eps>0$ was arbitrary.

Now, by the G\"{a}rtner--Ellis theorem, under
$\P\{\cdot| L^1=\omega_n\}$ the measure $L^2$ obeys a large
deviation principle on $\tilde{\skrim}_*(\skrix\times\skrix)$ with
rate function $\hat{I}_{\omega}(\varpi)=\sfrac12 \sup_{g\in\skric_2}\{
\langle
\varpi,g\rangle+\langle C\omega\otimes\omega, (1-e^{g})\rangle
\}$.
\end{pf}

Next, we express the rate function in term of relative entropies
and consequently show that it is
a good rate function. Recall the definition of the function
${I}_{\omega}$ from Theorem \ref{main}(a).
\begin{lemma}\label{randomg.Vrate}
\begin{longlist}
\item$\hat{I}_{\omega}(\varpi)=I_{\omega}(\varpi)$, for any
$\varpi\in\tilde{\skrim}_*(\skrix\times\skrix)$,
\item$I_{\omega}(\varpi)$ is a good rate function and
\item${\mathfrak{H}_C}(\varpi\| \omega)\ge0$ with
equality if and only if
$\varpi=C\omega\otimes\omega$.
\end{longlist}
\end{lemma}
\begin{pf}
(i) Suppose that $\varpi\not\ll C\omega\otimes\omega$.
Then, there exists
$a_0,b_0 \in\skrix$ with $C\omega\otimes\omega(a_0,b_0)=0$ and
$\varpi(a_0,b_0)>0$. Define $\hat{g}\dvtx\skrix\times\skrix
\rightarrow\R$ by
\[
\hat{g}(a,b)=\log\bigl[K\bigl(\one_{(a_0,b_0)}(a,b)+\one_{(b_0,a_0)}(a,b)\bigr)+1 \bigr]
\qquad\mbox{for $a,b\in\skrix$ and $K>0$.}
\]
For this choice of $\hat{g}$ we have
\[
\tfrac12\langle\varpi, \hat g \rangle+ \tfrac12 \langle C \omega
\otimes\omega, 1-e^{-\hat g} \rangle
\ge\tfrac12 \log(K+1)\varpi(a_{0},b_{0}) \to\infty\qquad \mbox{for
$K\uparrow\infty$.}
\]
Now suppose that $\varpi\ll C\omega\otimes\omega$. We have
\[
\hat{I}_{\omega}(\varpi) =
\sfrac{1}{2} \| C \omega\otimes\omega\| + \sfrac12
\sup_{g\in\skric_{2}} \{ \langle\varpi, g \rangle- \langle C \omega
\otimes\omega, e^g \rangle\}.
\]
By the substitution $h=e^{g} \frac{C \omega\otimes
\omega}{\varpi}$ the supremum equals
\begin{eqnarray*}
\mathop{\sup_{h\in\skric_{2}}}_{h \ge0} \biggl\langle
\varpi, \log\biggl( h \frac{\varpi}{C \omega\otimes\omega} \biggr)-h \biggr\rangle
& = &\mathop{\sup_{h\in\skric_{2}}}_{h \ge0} \langle\varpi, \log h -h
\rangle
+ \biggl\langle\varpi, \log\biggl(\frac{\varpi}{C \omega\otimes\omega} \biggr) \biggr\rangle\\
& = &- \| \varpi\| + H(\varpi \| C \omega\otimes\omega),
\end{eqnarray*}
where we have used $\sup_{x>0} \log x - x = -1$ in the last step.
This yields that $\hat{I}_{\omega}(\varpi)={I}_{\omega}(\varpi)$.

{\smallskipamount=0pt
\begin{longlist}[(iii)]
\item[(ii)] Recall that $I_{\omega}=\hat{I}_{\omega}$ is a rate function.
Moreover, for all $\alpha<\infty$, the level sets
$\{\varpi\in\tilde{\skrim}_*(\skrix\times\skrix)\dvtx
\sfrac12 {\mathfrak{H}_C}(\varpi\| \omega)\le\alpha\}$ are
bounded and therefore compact, so $I_{\omega}$ is a good rate function.

\item[(iii)] Consider the nonnegative function $\xi(x)=x\log x-x+1$, for
$x> 0$, $\xi(0)=1$, which has its only root at $x=1$. Note that
%
%
\begin{equation}\label{randomg.calculusnot}
{\mathfrak{H}_C}(\varpi\| \omega)=
\cases{\displaystyle\int\xi\circ g \,dC\omega\otimes\omega, &\quad if
$g:=\displaystyle\sfrac{d\varpi}{dC\omega\otimes\omega}\ge0$ exists, \cr
\infty, &\quad otherwise.}
\end{equation}
Hence ${\mathfrak{H}_C}(\varpi\| \omega)\ge0$, and, if
$\varpi=C\omega\otimes\omega$, then $\xi(\sfrac{d\varpi}{dC\omega\otimes
\omega})=\xi(1)=0$
and so ${\mathfrak{H}_C}(C\omega\otimes\omega\| \omega)=0$.
Conversely, if
${\mathfrak{H}_C}(\varpi\| \omega)=0$, then $\varpi(a,b)>0$ implies
$C\omega\otimes\omega(a,b)>0$, which then implies $\xi\circ g(a,b)=0$
and further $g(a,b)=1$. Hence $\varpi=C\omega\otimes\omega$, which
completes the proof of (iii).\qed
\end{longlist}}
\noqed\end{pf}

\section{\texorpdfstring{Proof of Theorem \protect\ref{randomg.LDprob} by the method of
types}{Proof of Theorem 2.5 by the method of
types}}\label{proofmain3}

Throughout the proof we may assume that $\omega(a)>0$ for all
$a\in\skrix$. It is easy to see that the law of the randomly
colored graph conditioned to have empirical color measure
$\omega_n$ and empirical pair measure $\varpi_n$
\[
\prob_{(\omega_n,\varpi_n)}:=\prob\{ \cdot | \Phi(M)=(\omega
_n,\varpi_n)\}
\]
can be described in the following manner:
\begin{itemize}
\item assign colors to the vertices by sampling without replacement
from the collection
of $n$ colors, which contains any color $a\in\skrix$ exactly $n\omega
_n(a)$ times;
\item for every unordered pair $\{a,b\}$ of colors create exactly
$n(a,b)$ edges by sampling
without replacement from the pool of possible edges connecting vertices
of color $a$ and $b$,
where
%
%
\begin{equation}\label{nabdef}
n(a,b):= \cases{n \varpi_n(a,b), &\quad if $a\not=b$,\vspace*{2pt}\cr
\dfrac n2 \varpi_n(a,b), &\quad if $a=b$.}
\end{equation}
\end{itemize}
We would like to reduce the calculation of probabilities to the
counting of
suitable configurations. To this end we introduce a numbering system,
which specifies, for each $\{a,b\}$, the
order in which edges are drawn in the second step. More precisely, the
edge-number $k$ is attached to both
vertices connecting the $k$th edge. Note that the total number of
edge-numbers attached to every
vertex corresponds to the degree of the vertex in the graph.
All permitted numberings are equally probable.

Denote by $Y^{\{a,b\}}_{^j}$ be the $j$th
edge drawn in the process of connecting vertices of colors $\{a,b\}$.
Let $\skria_n(\omega_n,\varpi_n)$
be the set of all possible configurations
\[
\bigl( \bigl(X(v) \dvtx v\in V \bigr);
\bigl(Y^{\{a,b\}}_k \dvtx k=1, \ldots, n(a,b);
\{a,b\}\subset\skrix\bigr) \bigr) ,
\]
and let
$\skrib_n(\omega_n,\varpi_n)$ be the set of all colored graphs $x$
with $L^1(x)=\omega_n$ and $L^2(x)=\varpi_n$. Define $\Psi\dvtx
\skria_n(\omega_n,\varpi_n)\to\skrib_n(\omega_n,\varpi_n)$ as the
canonical mapping which associates the colored graph to any
configuration, that is, ``forgets'' the numbering of the edges. Finally,
define
\[
\skrik^{\ssup n}(\omega_n,\varpi_n):= \{ {M}(x)
\mbox{ for some } x\in\skrib_n(\omega_n,\varpi_n) \}
\]
to be the set of all empirical neighborhood measures $M(x)$ arising
from colored
graphs $x$ with $n$ vertices with $\Phi(M(x))=(\omega_n, \varpi_n)$.
For any $\nu_n\in\skrik^{\ssup n}(\omega_n,\varpi_n)$ we have
%
%
\begin{equation}\label{randomg.probcom}\quad
{\prob} \{{M}=\nu_n | \Phi({M})
=(\omega_n,\varpi_n) \}=\frac{\sharp\{\tilde{x}\in\skria_n(\omega
_n,\varpi_n) \dvtx
{M \circ\Psi}(\tilde{x})=\nu_n \}}{\sharp\{\tilde{x}\in\skria_n(\omega
_n,\varpi_n) \}} .
\end{equation}
In our proofs we use the following form of Stirling's formula; see
\cite{Fel67}, page 54:
for all $n\in\N$,
\[
n^{n} e^{-n}\le n!\le(2\pi n)^{1/2} n^{n} e^{-n+1/(12n)} .
\]

\subsection{A bound on the number of empirical neighborhood
measures}\label{sec41}

In this section we provide an upper bound on the number of measures
in $\skrik^{\ssup n}(\omega_n,\varpi_n)$. We write $m$ for the number of
elements in $\skrix$.
\begin{lemma}\label{randomg.boundonK}
There exists $\vartheta>0$, depending on $m$ such that, if
$\omega_n\in\skrim_n(\skrix)$ and
$\varpi_n\in\tilde\skrim_{*,n}(\skrix\times\skrix)$, then
\[
\sharp\skrik^{\ssup n}(\omega_n,\varpi_n)\le\exp\bigl[\vartheta(\log
n) (n\|\varpi_n\|)^{({2m-1})/({2m})} \bigr].
\]
\end{lemma}

The proof is based on counting integer partitions of vectors.
To fix some notation, let ${\mathfrak I}_m={ (\N\cup\{0\} )}^{m}$ be the
collection of (nonnegative) integer vectors of length $m$. For any $\ell
\in{\mathfrak I}_m$
we denote by $\|\ell\|$ its magnitude, that is, the sum of its entries.

We introduce an ordering $\curlyeqsucc$ on ${\mathfrak I}_m$ such that,
for any
vectors
\[
\ell_{1}=\bigl(\ell_1^{ \ssup1},\ldots,\ell_1^{\ssup m}\bigr) \quad\mbox{and}\quad \ell
_{2}=\bigl(\ell_2^{\ssup1},\ldots,\ell_2^{\ssup m}\bigr),
\]
we write $\ell_1 \curlyeqsucc\ell_2$ if either:
\begin{longlist}
\item$\|\ell_{1}\|>\|\ell_{2}\|$, or
\item$\|\ell_{1}\|=\|\ell_{2}\|$ and there is $j\in\{1,\ldots
,m\}$ with
$\ell_{1}^{\ssup k}=\ell_{2}^{\ssup k}$, for all $k<j$, and
$\ell_{1}^{\ssup j}>\ell_{2}^{\ssup j}$ or
\item$\ell_1=\ell_2$.
\end{longlist}
A collection $(\ell_{1}, \ldots, \ell_{k})$ of elements in
${\mathfrak I}_m$ is an \textit{integer partition} of the vector $\ell
\in
{\mathfrak I}_m$, if
\[
\ell_{1}\curlyeqsucc\cdots\curlyeqsucc\ell_{k} \not= 0
\quad\mbox{and}\quad \ell_{1}+\cdots+\ell_{k}=\ell.
\]
Any integer partition of a vector $\ell\in{\mathfrak I}_m$ induces an integer
partition $\|\ell_{1}\|,\ldots,\break\|\ell_{k}\|$ of its magnitude $\|\ell\|
$, which
we call its \textit{sum-partition}. We denote by $\skrip_m(\ell)$ the
set of
integer partitions of $\ell$.
\begin{lemma}\label{randomg.dpartition}
There exists $\vartheta>0$, which depends on $m$ such that, for any
$\ell\in{\mathfrak I}_m$ of magnitude $n$,
\[
\sharp\skrip_m(\ell) \le\exp\bigl[\vartheta(\log n) n^{({2m-1})/({2m})}
\bigr] .
\]
\end{lemma}
\begin{pf}
Let $\ell\in{\mathfrak I}_m$ be a vector of magnitude $n$ and
$(\ell_1, \ldots, \ell_k)$ be an integer partition of $\ell$. We
relabel the partition as $(\emm_{1,1}, \ldots, \emm_{1,k_1};
\emm_{2,1}, \ldots, \emm_{2,k_2};\break \ldots; \emm_{r,1},\ldots,
\emm_{r,k_r})$ such that all vectors in the same block (indicated by
the first subscript) have the same magnitude, which we denote
$y_1,\ldots, y_r$, and such that $y_1>\cdots>y_r>0$. Note that for
the block sizes we have $k_1+\cdots+k_r=k$ and that $k_1y_1+\cdots
+k_ry_r=n$.

For a moment, look at a fixed block $\emm_{j,1}, \ldots, \emm_{j,k_j}$.
It is easy to see that the number
of integer vectors of length $m$ and magnitude $y_j$ is given by
\[
b(y_j,m):=\pmatrix{y_j+m-1\cr m-1} \le c(m) y_j^{m-1} .
\]
Writing $\emm_{j,0}$ for the largest and $\emm_{j,k_j+1}$ for the smallest
of these vectors in the ordering $\curlyeqsucc$, we note that
\[
p\dvtx\{0,\ldots, k_j+1\} \to\{ \emm\in{\mathfrak I}_m \dvtx\|
\emm\|=y_j \},\qquad
p(i)=\emm_{j,i},
\]
is a nonincreasing path of length $k_j+2$ into an ordered set of size
$b(y_j,m)$,
which connects the smallest to the largest element. The number of such
paths is easily
seen to be
\[
\pmatrix{b(y_j,m)+k_j\cr k_j}.
\]

Therefore, the number of integer partitions of $\ell$ with given
sum-partition $(y_1, \stackrel{k_1}{\ldots}\,,y_1,\ldots,y_r,\stackrel
{k_r}{\ldots}\,, y_r)$ is
\[
\prod_{j=1}^{r}\pmatrix{b(y_j,m)+k_{j}\cr k_j}
\le\mathop{\max_{a_1,\ldots, a_r>0}}_{\sum a_j=n} \prod_{j=1}^{r} \biggl\{
\mathop{\max_{y,k \in\N}}_{yk=a_j}
\pmatrix{c(m)y^{m-1}+k\cr k} \biggr\}.
\]

To maximize the binomial coefficient over the set $yk=a_j$, we
distinguish between the cases when (i) $a_j\le c(m) y^{m}$,
(ii) $a_j>c(m)y^{m}$ and observe that
\[
\pmatrix{c(m)y^{m-1}+\dfrac{a_j}{y}\cr\dfrac{a_j}{y}}\le
\cases{\pmatrix{2c(m)y^{m-1}\vspace*{2pt}\cr\dfrac{a_j}{y}}, &\quad if $a_{j}\le c(m) y^{m}$,
\vspace*{6pt}\cr
\pmatrix{2\dfrac{a_j}{y}\vspace*{2pt}\cr c(m)y^{m-1}}, &\quad if
$a_{j}>c(m)y^{m}$.}
\]

\textit{Case} (i): Using the upper bound ${i\choose r}\le(\sfrac{ie}{r}
)^{r}$, for $r,i\in\N$
with $r\le i $ and the inequality $
(\sfrac{a_j}{c(m)})^{1/m}\le y\le a_j\le n$ we obtain, for some
constants $C_0=C_0(m)$, $C_1=C_1(m)$,
\begin{eqnarray*}
\pmatrix{2c(m)y^{m-1}\cr a_{j}/y} &\le&
\biggl(\frac{2c(m)y^{m-1}e}{a_{j}/y} \biggr)^{a_{j}/y}\le
\exp\bigl( C_0 (a_{j}/y) \log n \bigr)\\
&\le&
\exp\bigl((\log n) C_1 a_{j}^{({m-1})/{m}} \bigr) .
\end{eqnarray*}

\textit{Case} (ii): The same upper bound for binomial
coefficients and $ 1\le y\le(\sfrac{a_j}{c})^{1/m}\le
a_j\le n$ yield for some constant $C_2=C_2(m)$,
\[
\pmatrix{2a_{j}/y\cr c(m)y^{m-1}}\le
\biggl(\frac{2(a_{j}/y)e}{c(m)y^{m-1}} \biggr)^{c(m)y^{m-1}}\le
\exp\bigl((\log n) C_2 a_{j}^{({m-1})/{m}} \bigr) .
\]

From these cases, we have for some $C=C(m)>0$, the upper bound
\[
\prod_{j=1}^{r}\pmatrix{b(y_j,m)+k_{j}\cr k_j}
\le\mathop{\max_{a_1,\ldots, a_r>0}}_{\sum a_j=n} \prod_{j=1}^r
\exp\bigl((\log n) C a_j^{({m-1})/{m}} \bigr) ,
\]
which is estimated
further (using H\"older's inequality) by
\[
\exp\Biggl( (\log n) C \sum_{j=1}^ra_j^{({m-1})/{m}} \Biggr)
\le\exp\Biggl( (\log n) C r^{1/{m}} \Biggl( \sum_{j=1}^ra_j \Biggr)^{
({m-1})/{m}} \Biggr).
\]
We observe that all $y_j$ are different, positive and that their sum is
not greater than~$n$, so
we have that
\[
r^2/2\le1+\cdots+r\le y_{1}+\cdots+y_{r}\le n.
\]
Recalling that
$a_1+\cdots+ a_r=n$, our upper bound becomes
\[
\exp\bigl( (\log n) C (2n)^{{1}/({2m})} n^{({m-1})/{m}} \bigr)
= \exp\biggl( \sfrac{\vartheta}{2} (\log n)
n^{({2m-1})/({2m})} \biggr)
\]
for some $\vartheta=\vartheta(m)>0$.
Note that from our argument so far one can easily recover the
well-known fact that the number of integer partitions of $n$ is
bounded by $e^{(\vartheta/2) \sqrt{n}}$, the full asymptotics being discovered
by Hardy and Ramanujan in 1918. Combining this with the
upper bound for the number of integer partitions with a given
sum-partition, we obtain the claim.
\end{pf}
\begin{pf*}{Proof of Lemma \ref{randomg.boundonK}} Suppose $\omega_n\in
\skrim_n(\skrix)$ and
$\varpi_n\in\tilde\skrim_{*,n}(\skrix\times\skrix)$. For $a\in\skrix$,
we look at the mappings
\[
\Phi_a\dvtx\skrik^{\ssup n}(\omega_n,\varpi_n)\ni\frac{1}{n}\sum
_{v\in V}\delta_{(X(v),L(v))}
\mapsto\bigl( L^a_1, \ldots, L^a_{n\omega(a)} \bigr) ,
\]
where $(L^a_1, \ldots, L^a_{n\omega(a)})$ is the ordering of the
vectors $L(v)$, for all
$v\in V$ with $X(v)=a$, and thus constitutes an integer partition of
the vector
$( n\varpi_n(a,b) \dvtx b\in\skrix)$,
which has magnitude $n\sum_b \varpi_n(a,b)$. The combined mapping $\Phi
=(\Phi_a \dvtx a\in\skrix)$
is injective, and therefore, by Lemma \ref{randomg.dpartition},
\begin{eqnarray*}
\sharp\skrik^{\ssup n}(\omega_n,\varpi_n)& \le &\exp\biggl[\vartheta
\sum_{a\in\skrix} \log\biggl( n\sum_{b\in\skrix} \varpi_n(a,b) \biggr)
\biggl( n\sum_{b\in\skrix} \varpi_n(a,b) \biggr)^{({2m-1})/({2m})} \biggr]
\\
&\le&\exp\bigl[\vartheta m \log(n\|\varpi_n\| )
(n\|\varpi_n\| )^{({2m-1})/({2m})} \bigr]
,
\end{eqnarray*}
where we have used the fact that $\sum_b\varpi_n(a,b)\le\|\varpi_n\|$
in the last step.
\end{pf*}

\subsection{\texorpdfstring{Proof of the upper bound in Theorem \protect\ref
{main3}}{Proof of the upper bound in Theorem 2.5}}\label{sec42}

We are now ready to prove an upper bound for the large deviation probability
in Theorem \ref{main3}. 
%
\begin{lemma}\label{upperb}
For any sequence $(\nu_n)$ with $\nu_n\in\skrik^{\ssup n}(\omega
_n,\varpi_n)$ we have
\[
{\prob} \{M=\nu_n | \Phi(M)=(\omega_n,\varpi_n) \}
\le\exp\bigl( -nH(\nu_n \| Q_n)+\eps_{1}^{\ssup n}(\nu_n) \bigr) ,
\]
where
\[
Q_n(a , \ell)=\omega_n(a)\prod_{b\in\skrix}
\frac{e^{-\varpi_n(a,b)/\omega_n(a)}[\varpi_n(a,b)/\omega_n(a)]^{\ell
(b)}}{\ell(b)!}\qquad
\mbox{for $\ell\in\skrin(\skrix)$}
\]
and
\[
\lim_{n\uparrow\infty} \frac1n \sup_{\nu_n\in\skrik^{\ssup n}(\omega
_n, \varpi_n)}
\eps_{1}^{\ssup n}(\nu_n)=0 .
\]
\end{lemma}
\begin{pf}
The proof of this lemma is based on the method of types; see \cite
{DZ98}, Chapter 2.
Recall from \eqref{randomg.probcom} that, for any $\nu_n\in\skrik
^{\ssup n}(\omega_n,\varpi_n)$, we have
\[
{\prob} \{{M}=\nu_n | \Phi({M})
=(\omega_n,\varpi_n) \}=\frac{\sharp\{\tilde{x}\in\skria_n(\omega
_n,\varpi_n) \dvtx
{M \circ\Psi}(\tilde{x})=\nu_n \}}{\sharp\{\tilde{x}\in\skria_n(\omega
_n,\varpi_n) \}} .
\]
Now, by elementary counting, the denominator on the right-hand side of
\eqref
{randomg.probcom} is
%
%
\begin{equation}\label{denom}\qquad
\pmatrix{n\cr n\omega_n(a), a\in\skrix}
\prod_{\{a,b\}} \prod_{k=1}^{n(a,b)}
\biggl( \frac{n^2\omega_n(a)\omega_n(b) - n\omega_n(a)\one_{\{a=b\}}}{1+\one
_{\{a=b\}}} - (k-1) \biggr) .
\end{equation}
For a given empirical neighborhood measure $\nu_n\in\skrik^{\ssup
n}(\omega_n,\varpi_n)$
the numerator is probably too tricky to find explicitly. However, an
easy upper bound is
%
%
\begin{eqnarray}\label{numerupper}
&&\pmatrix{n\cr n\nu_n(a,\ell), a\in\skrix, \ell\in\skrin(\skrix)}
2^{-{n}/2 \varpi_n(\Delta)}\nonumber\\[-8pt]\\[-8pt]
&&\qquad{}\times
\prod_{(a,b)} \pmatrix{n \varpi_n(a,b)\cr\ell_{a}^{\ssup j}(b),
j=1,\ldots,n\omega_n(a)},\nonumber
\end{eqnarray}
where $\ell_{a}^{\ssup j}(b), j=1,\ldots,n\omega_n(a)$ are any
enumeration of the family
containing each $\ell(b)$ with multiplicity $n\nu_n(a,\ell)$. This can
be seen from the following
construction: first allocate to each vertex some $(a,\ell)\in\skrix
\times\skrin(\skrix)$ in such a
way that every vector $(a,\ell)$ is allocated $n\nu_n(a,\ell)$ times.
The first binomial coefficient
in \eqref{numerupper} represents the number of possible ways to do
this. For any $(a,b)\in\skrix\times\skrix$
distribute the numbers in $\{1,\ldots,n \varpi_n(a,b)\}$ among the
vertices with color $a$ so that
a vertex carrying vector $(a,\ell)$ gets exactly $\ell(b)$ numbers.
Once this is done for both $(a,b)$
and $(b,a)$, each vertex of color $a$ or $b$ carries a set of numbers;
if $a\not=b$ each number in
$\{1,\ldots, n \varpi_n(a,b)\}$ occurs exactly twice in total, if $a=b$
it occurs exactly once. Next, for $k=1,\ldots, n(a,b)$,
if $a\not=b$ draw the $k$th edge between the two vertices of color $a$
and $b$
carrying number $k$, if $a=b$ draw the $k$th edge between the vertices
with number $k$ and $2k$. The remaining
factor in \eqref{numerupper} represents the number of possible ways to
do this, with the power of two
discounting the fact that for edges connecting vertices of the same
color two numbering schemes lead to
the same configuration. By this construction, every element $\tilde
{x}\in\skria_n(\omega_n,\varpi_n)$ with $M \circ\Psi(\tilde{x})=\nu_n$
has been constructed exactly once, but also some graphs with loops or
multiple edges can occur,
so that \eqref{numerupper} is an upper bound for the numerator in \eqref
{randomg.probcom}.

Combining \eqref{randomg.probcom}, \eqref{denom}
and \eqref{numerupper} we get
%
%
\begin{eqnarray}\label{Tclass}
&&{\prob} \{{M}=\nu_n | \Phi({M})=(\omega_n,\varpi_n) \}\nonumber\\
&&\qquad\le\prod_{a\in\skrix}\pmatrix{n\omega_n(a)\cr n\nu_n(a,\ell), \ell\in
\skrin(\skrix)}\nonumber\\[-8pt]\\[-8pt]
&&\qquad\quad{}\times
\prod_{(a,b)}\pmatrix{n\varpi_n(a,b)\cr\ell_{a}^{\ssup j}(b),
j=1,\ldots,n\omega_n(a)}2^{-n/2 \varpi_n(\Delta)}\nonumber\\
&&\qquad\quad{} \times \prod_{\{a,b\}}
\prod_{k=1}^{n(a,b)}
\biggl( \sfrac{n^2\omega_n(a)\omega_n(b) - n\omega_n(a)\one_{\{a=b\}}}{1+\one
_{\{a=b\}}} - (k-1) \biggr)^{-1}
.\nonumber
\end{eqnarray}

It remains to analyze the asymptotics of this upper bound. Using
Stirling's formula, we obtain
\begin{eqnarray*}
&&\prod_{a\in\skrix}\pmatrix{n\omega_n(a)\cr n\nu_n(a,\ell), \ell\in
\skrin(\skrix)}
\\
&&\qquad \le\exp\biggl(n\sum_{a}\omega_n(a)\log\omega_n(a)-n\sum_{(a,\ell)}\nu
_n(a,\ell)
\log\nu_n(a,\ell) \biggr)\\
&&\qquad\quad{} \times\exp\biggl( \sfrac{m}{2}\log(2\pi n)+\sfrac{1}{n}
\sum_{a}\sfrac{1}{12\omega_n(a)} \biggr).
\end{eqnarray*}
We observe that
\[
\prod_{j=1}^{n\omega_n(a)} \bigl(\ell_{a}^{\ssup j}(b) \bigr)!=
\exp\biggl(n\sum_{\ell}\log(\ell(b)! ) \nu_n(a, \ell) \biggr),
\]
and hence
\begin{eqnarray*}
\hspace*{-3pt}&&\pmatrix{n\varpi_n(a,b)\cr\ell_{a}^{\ssup j}(b), j\le n\omega_n(a)}
\\
\hspace*{-3pt}&&\qquad \le\exp\biggl(-n\sum_{\ell} \log(\ell(b)! ) \nu_n(a,\ell)
+n\varpi_n(a,b)\log(n\varpi_n(a,b) )-n\varpi_n(a,b) \biggr) \\
\hspace*{-3pt}&&\qquad\quad{} \times\exp\biggl(\sfrac{1}{12n\varpi_n(a,b)}+\sfrac
{1}{2}\log(2\pi n) \biggr).
\end{eqnarray*}
Next, we obtain
\begin{eqnarray*}
&&\prod_{k=1}^{n(a,b)} \biggl( \sfrac{n^2\omega_n(a)\omega_n(b) - n\omega
_n(a)\one_{\{a=b\}}}{1+\one_{\{a=b\}}} - (k-1) \biggr)
\\
&&\qquad \ge\exp\biggl( n(a,b) \log\biggl(\sfrac{n^2\omega_n(a)\omega_n(b)}{1+\one
_{\{a=b\}}} \biggr) \biggr) \\
&&\qquad\quad{} \times\exp\biggl( n(a,b)
\log\biggl(1- \sfrac{\one\{a=b\}}{2n\omega_n(a)}- \sfrac{2n(a,b)}{n^2\omega
_n(a)\omega_n(b)} \biggr) \biggr) .
\end{eqnarray*}
Putting everything together and denoting by $H(\omega)=-\sum_{y\in\skriy
} \omega(y) \log\omega(y)$
the entropy of a measure $\omega\in\skrim(\skriy)$, we get
\begin{eqnarray*}
&&{\prob} \{M=\nu_n | \Phi(M)=(\omega_n,\varpi_n) \}\\
&&\qquad\le\exp\biggl(-nH(\omega_n)+nH(\nu_n)-n\sum_{(a,b)}\sum_{\ell}(\log\ell
(b)!) \nu_n(a,\ell)\\
&&\qquad\hspace*{31.7pt}{}
+ n\sum_{(a,b)}\varpi_n(a,b) \log\varpi_n(a,b) -n\sum_{(a,b)}\varpi_n(a,b)\\
&&\qquad\hspace*{31.7pt}{} - \sfrac n2 \sum_{(a,b)}\varpi
_n(a,b) \log\biggl(
\sfrac{\omega_n(a)\omega_n(b)}{1+\one_{\{a=b\}}} \biggr)-\sfrac
n2 \varpi_n(\Delta)\log2+\varepsilon_{1}^{\ssup n} \biggr) ,
\end{eqnarray*}
for a sequence $\eps_1^{\ssup n}$ which does not depend on $\nu_n$ and satisfies
$\lim_{n\uparrow\infty} \sfrac1n \eps_{1}^{\ssup n}=0$. To give
the right-hand side the form
as stated in the theorem, we observe that
\begin{eqnarray*}
\hspace*{-5pt}&&
H(\omega_n)-H(\nu_n) +\sum_{(a,b)}\sum_{\ell}(\log\ell(b)!) \nu
_n(a,\ell)
- \sum_{(a,b)}\varpi_n(a,b) \log\varpi_n(a,b)\\[-1.2pt]
\hspace*{-5pt}&&\quad{} +\sum_{(a,b)}\varpi
_n(a,b) + \sfrac12 \sum_{(a,b) }\varpi_n(a,b) \log(\omega
_n(a)\omega_n(b) )\\[-1.2pt]
\hspace*{-5pt}&&\qquad=\sum_{(a,\ell)}\nu_n(a,\ell) \biggl[\log\nu_n(a,\ell)-\log\omega_n(a)\\[-1.2pt]
\hspace*{-5pt}&&\qquad\hspace*{68.2pt}{}-\sum
_{b} \biggl(\log\biggl(\sfrac{\varpi_n(a,b)}
{\omega_n(a)}\biggr)^{\ell(b)}-\sfrac{\varpi_n(a,b)}{\omega_n(a)}-(\log\ell
(b)!) \biggr) \biggr]\\
\hspace*{-5pt}&&\qquad=\sum_{(a,\ell)}\nu_n(a,\ell) \biggl[\log\nu_n(a,\ell)\\
\hspace*{-5pt}&&\qquad\quad\hspace*{59.3pt}-\log\biggl(\omega_n(a)\\
\hspace*{-5pt}&&\qquad\quad\hspace*{87.7pt}{}\times\prod_{b}
\sfrac{(\varpi_n(a,b)/\omega_n(a))^{\ell(b)}\exp(-\varpi_n(a,b)/\omega
_n(a))}{\ell(b)!} \biggr) \biggr]\\
\hspace*{-5pt}&&\qquad=H(\nu_n \| Q_n),
\end{eqnarray*}
which completes the proof of Lemma \ref{upperb}.
\end{pf}

We can now complete the proof of the upper bound in Theorem \ref{main3}
by combining
Lemmas \ref{randomg.boundonK} and \ref{upperb}. Suppose that
$\Gamma\subset\skrim(\skrix\times\skrin(\skrix))$ is a closed set. Then
\begin{eqnarray*}
&&{\prob} \{M \in\Gamma | \Phi(M)=
(\omega_n,\varpi_n) \}\\
&&\qquad=\sum_{\nu_n\in\Gamma\cap\skrik^{\ssup n}(\omega
_n,\varpi_n)}{\prob} \{M=\nu_n |
\Phi(M)=(\omega_n,\varpi_n) \}\\
&&\qquad\le\sharp\skrik^{\ssup n}(\omega_n,\varpi_n) \exp\Bigl( -n \inf_{\nu_n
\in\Gamma\cap\skrik^{\ssup n}(\omega_n,\varpi_n)}
H(\nu_n \| Q_n)\\
&&\qquad\quad\hspace*{115.2pt}{} + \sup_{\nu_n\in\skrik^{\ssup n}(\omega_n,\varpi_n)}
\eps^{\ssup n}_1(\nu_n) \Bigr) .
\end{eqnarray*}
We have already seen that $\sfrac1n \sup_{\nu_n} \eps^{\ssup
n}_1(\nu_n)$ and $\sfrac1n
\log\sharp\skrik^{\ssup n}(\omega_n,\varpi_n)$ converge to zero. It
remains to check that
%
%
\begin{equation}\label{equ.entropyproof}
{\lim_{n\to\infty} \sup_{\nu_n\in\skrik^{\ssup n}(\omega_n,\varpi_n)}}
| H(\nu_n \| Q_n) - H(\nu_n \| Q) | = 0 .
\end{equation}
To do this, we observe that
%
%
\begin{eqnarray}\label{randomg.entt}
&&H(\nu_n \| Q_n) -
H(\nu_n \| Q)
\nonumber\\
&&\qquad=\sum_{(a,\ell)\in\skrix\times\skrin(\skrix)}\nu_n(a,\ell)\log\sfrac
{Q(a,\ell)}{Q_n(a,\ell)}\nonumber\\[-8pt]\\[-8pt]
&&\qquad=-H(\omega_n \| \omega)-H(\varpi_n \| \varpi)-\sum_{a,b\in\skrix
}\varpi(a,b)\sfrac{\omega_n(a)}{\omega(a)}
\nonumber\\
&&\qquad\quad{}+\sum_{a,b\in\skrix}\varpi(a,b)\log
\sfrac{\omega_n(a)}{\omega(a)}+\|\varpi_n\|.\nonumber
\end{eqnarray}
Note that this expression does not depend on $\nu_n$.
As the first, second and fourth term of \eqref{randomg.entt} converge
to $0$, and the third
and fifth term converge to $\|\varpi\|$, the expression \eqref{randomg.entt}
vanishes in the limit and this completes the proof of the upper bound
in Theorem \ref{main3}.

\subsection{An upper bound on the support of empirical neighborhood
measures}\label{sec43}

The cardinality of the support, denoted $\sharp\skris(\nu)$, of an
empirical neighborhood measure $\nu$ of a graph with $n$ vertices
is naturally bounded by $n$. For the proof of the lower bound in
Theorem \ref{main3} we need a better upper bound. We still use $m$
to denote the cardinality of $\skrix$, and let
\[
C:=2^m \frac{\Gamma(m+2)^{{m}/({m+1})}}{\Gamma(m)}
\quad\mbox{and}\quad D:=2^m \frac{(m+1)^{m}}{\Gamma(m)},
\]
where $\Gamma(\cdot)$ is the Gamma function.
\begin{lemma}\label{randomg.finiteness}
For every $(\omega_n,\varpi_n)\in\skrim_n(\skrix)\times\tilde\skrim
_{*,n}(\skrix\times\skrix)$
and $\nu_n\in\skrim_n(\skrix\times\skrin(\skrix))$ with $\Phi(\nu
_n)=(\omega_n,\varpi_n)$, we have
%
%
\begin{equation}\label{randomg.inequality1}
\sharp\skris(\nu_n)\le C [n\|\varpi_n\| ]^{{m}/({m+1})} + D.
\end{equation}
\end{lemma}

The following lemma provides a step in the proof of Lemma \ref
{randomg.finiteness}.
\begin{lemma}\label{randomg.finiteness1}
Suppose $j\in\N\cup\{0\}$ and $n\in\N$. Then,
%
%
\begin{eqnarray}\label{randomg.eqnumber}
\sfrac{1}{\Gamma(n)} j^{n-1}&\le&
\sharp\bigl\{(l_1, \ldots, l_n)\in(\N\cup\{0\} )^{n}\dvtx
l_1+\cdots+l_n=j \bigr\}\nonumber\\[-8pt]\\[-8pt]
&\le&\sfrac{1}{\Gamma(n)} (j+n)^{n-1}.\nonumber
\end{eqnarray}
\end{lemma}
\begin{pf}
The proof is by induction on $n$. Equation \eqref{randomg.eqnumber}
holds trivially
for all $j\in\N\cup\{0\}$ and $n=1,2$, so we assume it holds for all
$j$ and
$n\ge2$. By the induction hypothesis, for any $j$,
\begin{eqnarray*}
&&\frac{1}{\Gamma(n)}\sum_{l=0}^{j}(j-l)^{n-1}\\
&&\qquad \le \sum_{l=0}^{j} \sharp\bigl\{(l_1, \ldots, l_{n-1})\in(\N\cup\{0\}
)^{n-1}\dvtx
l_1+\cdots+l_{n-1}=j-l \bigr\} \\
&&\qquad = \sharp\bigl\{(l_1, \ldots, l_n)\in(\N\cup\{0\} )^{n}\dvtx
l_1+\cdots+l_n=j \bigr\}\\
&&\qquad\le\frac{1}{\Gamma(n)} \sum_{l=0}^{j}(j-l+n)^{n-1}.
\end{eqnarray*}
For the first and last term, we obtain the lower and upper
bounds
\[
\sum_{l=0}^{j}(j-l)^{n-1}\ge
\int_{0}^{j}y^{n-1} \,dy=\frac{1}{n} j^{n}=\frac{\Gamma(n)}{\Gamma(n+1)}j^{n}
\]
and
\begin{eqnarray*}
\sum_{l=0}^{j}(j-l+n)^{n-1}&\le&\int_{n}^{j+n}y^{n-1}\,dy\le
\int_{0}^{j+n+1}y^{n-1}\,dy=\frac{1}{n}
(j+n+1)^{n}\\
&=&\frac{\Gamma(n)}{\Gamma(n+1)}(j+n+1)^{n},
\end{eqnarray*}
which
yields inequality \eqref{randomg.eqnumber} for $n+1$ instead of $n$, and
completes the induction.
\end{pf}
\begin{pf*}{Proof of Lemma \ref{randomg.finiteness}}
Suppose $(\omega_n,\varpi_n)\in\skrim_n(\skrix)\times\skrim_n(\skrix
\times\skrix)$. Let
\begin{eqnarray*}
a_m(j):\!&=&\sharp\biggl\{ (a,\ell)\in\skrix\times\skrin(\skrix) \dvtx
\sum_{b\in\skrix} \ell(b)=j \biggr\}\\
&=&m\times\sharp\bigl\{(l_1, \ldots, l_m)\in(\N\cup\{0\} )^{m}\dvtx
l_1+\cdots+l_m=j \bigr\}.
\end{eqnarray*}
For any positive integer $k$ we write
\[
\theta_k=\min\Biggl\{\theta\in\N\dvtx\sum_{j=0}^{\theta}a_m(j)\ge k \Biggr\}.
\]
We observe from Lemma \ref{randomg.finiteness1} that
\begin{eqnarray*}
k&\le&\sum_{j=0}^{\theta_{k}}a_m(j)\le
m\sum_{j=0}^{\theta_{k}}\sfrac{1}{\Gamma(m)} (j+m)^{m-1}
\le\sfrac{m}{\Gamma(m)} \int_{0}^{\theta_{k}+m}y^{m-1}\,dy\\
&=&
\sfrac{1}{\Gamma(m)} (\theta_{k}+m )^{m}.
\end{eqnarray*}
Thus, we have $ \theta_{k}\ge
(k\Gamma(m) )^{1/m}
-m=:\alpha_k$.
This yields
%
%
\begin{eqnarray}\label{randomg.equ1}
\sum_{j=0}^{\theta_{k}}ja_m(j)&\ge&
\sfrac{1}{\Gamma(m)}\sum_{j=0}^{\lceil\alpha_{k}\rceil}j^{m}
\ge\sfrac{1}{\Gamma(m)}\int_{0}^{\alpha_{k}-1}y^{m}\,dy\nonumber\\[-8pt]\\[-8pt]
&\ge&\sfrac{1}{\Gamma(m+2)}
(\alpha_{k}-1 )^{m+1} ,\nonumber
\end{eqnarray}
where $\lceil y\rceil$ is the smallest integer greater or equal to $y$.

Observe that the size of the support of the measure $\nu_n\in\skrik
^{\ssup n}(\omega_n,\varpi_n)$
satisfies
\[
\sharp\skris(\nu_n) \le\max\Biggl\{
k\dvtx\sum_{j=0}^{\theta_{k}}ja_{m}(j)\le n\|\varpi_n\| \Biggr\},
\]
and hence, using \eqref{randomg.equ1} and the inequality $(a+b)^m\le
2^m(a^m + b^m)$ for $a,b\ge0$,
\begin{eqnarray*}
\sharp\skris(\nu_n)&\le&\max\biggl\{k\dvtx\sfrac{1}{\Gamma(m+2)}
(\alpha_{k}-1 )^{m+1}\le n\|\varpi_n\| \biggr\}
\\
& \le &\Gamma(m)^{-1} \bigl( (n\|\varpi_n\|)^{1/({m+1})} \Gamma
(m+2)^{1/({m+1})} + m + 1 \bigr)^m\\
&\le& C (n\|\varpi_n\|)^{{m}/({m+1})} + D ,
\end{eqnarray*}
where the constants $C$, $D$ were defined before the formulation of
the lemma.
\end{pf*}

\subsection{Approximation by empirical neighborhood measures}\label{sec44}

Throughout this section we assume that $\omega_n\in\skrim_n(\skrix)$
with $\omega_n\to\omega$,
$\varpi_n\in\tilde\skrim_{*,n}(\skrix\times\skrix)$ with
$\varpi_n\to\varpi$, and that
$(\varpi,{\nu})$ is sub-consistent and ${\nu}_1=\omega$.
Our aim is to show that $\nu$ can be approximated in the weak
topology by some $\nu_n\in\skrim_n(\skrix\times\skrin(\skrix))$
with $\Phi(\nu_n)=(\omega_n, \varpi_n)$
and the additional feature that
%
%
\begin{equation}\label{degr}
\sum_{b\in\skrix} \ell(b)\le n^{1/3} \qquad\mbox{for $\nu_n$-almost
every $(a,\ell)$.}
\end{equation}
The approximation will be done in three steps, given as Lemmas \ref
{newlemma}, \ref{randomg.approx}
and \ref{ntoathird}. We denote by $d$ the metric of total variation,
that is,
\[
d(\nu,\tilde{\nu})=\sfrac{1}{2}\sum_{(a,\ell)\in\skrix\times\skrin
(\skrix)}
|\nu(a,\ell)-\tilde{\nu}(a,\ell)|\qquad
\mbox{for }\nu,\tilde{\nu}\in\skrim\bigl(\skrix\times\skrin(\skrix)\bigr).
\]
This metric generates the weak topology.
\begin{lemma}[(Approximation step 1)]\label{newlemma}
For every $\eps>0$, there exist $\hat\nu\in\skrim(\skrix\times\skrin
(\skrix))$
and $\hat\varpi\in\tilde\skrim(\skrix\times\skrix)$
such that $|\varpi(a,b)-\hat\varpi(a,b)|\le\eps$ for all $a,b\in\skrix$,
$d(\nu,\hat\nu)\le\eps$ and $(\hat\varpi, \hat\nu)$ is consistent.
\end{lemma}
\begin{pf}
By our assumption $(\varpi,\nu)$ is sub-consistent. For any $b\in\skrix
$ define
$e^{\ssup b}\in\skrin(\skrix)$ by $e^{\ssup b}(a)=0$ if $a\not=b$, and
$e^{\ssup b}(b)=1$. For large $n$
define measures $\hat\nu_n\in\skrim(\skrix\times\skrin(\skrix))$ by
\begin{eqnarray*}
\hat\nu_n(a,\ell) &=& \nu(a,\ell)
\biggl(1- \sfrac{\|\varpi\|- \|\langle\nu(\cdot,\ell), \ell(\cdot)\rangle\|
}{n} \biggr)\\
&&{} + \sum_{b\in\skrix} \one\bigl\{ \ell=ne^{\ssup b} \bigr\}
\sfrac{\varpi(a,b)- \langle\nu(\cdot,\ell), \ell(\cdot)\rangle
(a,b)}{n} .
\end{eqnarray*}
Note that $\hat\nu_n \to\nu$ and that, for all $a,b\in\skrix$,
\begin{eqnarray*}
&&\sum_{\ell\in\skrin(\skrix)} \hat\nu_n(a,\ell)\ell(b)\\
&&\qquad= \biggl(1- \sfrac{\|\varpi\|- \|\langle\nu(\cdot,\ell), \ell(\cdot)\rangle\|
}{n} \biggr)
\sum_{\ell\in\skrin(\skrix)}\nu(a,\ell)\ell(b)\\
&&\qquad\quad{}+ \varpi(a,b)- \langle\nu(\cdot,\ell), \ell(\cdot)\rangle(a,b)\\
&&\qquad = \varpi(a,b) - \sfrac{\|\varpi\|- \|\langle\nu(\cdot,\ell), \ell
(\cdot)\rangle\|}{n}
\langle\nu(\cdot,\ell), \ell(\cdot)\rangle(a,b)\\
&&\qquad\stackrel{n\uparrow\infty}{\longrightarrow} \varpi(a,b) .
\end{eqnarray*}
Hence, defining $\hat\varpi_n$ by $\hat\varpi_n(a,b)=\sum\hat\nu
_n(a,\ell)\ell(b)$,
we have a sequence of consistent pairs $(\hat\varpi_n, \hat\nu_n)$
converging to
$(\varpi,\nu)$, as required.
\end{pf}
\begin{lemma}[(Approximation step 2)]\label{randomg.approx}
For every $\eps>0$, there
exists $n(\eps)$ such that, for all $n\ge n(\eps)$ there exists
$\nu_n\in\skrim_n(\skrix\times\skrin(\skrix))$ with
$\Phi(\nu_n)=(\omega_n,\varpi_n)$ such that $d(\nu_n,\nu)\le\eps$.
\end{lemma}

The key to the construction of the measure $\nu_n$ is the following
``law of large numbers.''
\begin{lemma}\label{random.approx1}
For every $\delta>0$, there exists $\hat{\nu}\in\skrim(\skrix\times
\skrin(\skrix))$ with
$d(\nu, \hat\nu)<\delta$ such that, for i.i.d. $\skrin(\skrix)$-valued
random variables ${\ell}^{a}_j, j=1,\ldots,n\omega_n(a)$ with law
$\hat{\nu}( \cdot | a) := {\hat{\nu}( \{a\} \times\cdot
)}/{\hat{\nu}_1(a)}$,
almost surely,
%
%
\begin{equation}\label{randomg.estovershoot1}
\limsup_{n\rightarrow\infty} \Biggl( \frac1n \sum_{j=1}^{n\omega_n(a)}
{\ell}_{j}^{a}(b)-\varpi_n(a,b) \Biggr)\le\delta\qquad\mbox{for all
$a,b\in\skrix$.}
\end{equation}
\end{lemma}
\begin{pf}
By Lemma \ref{newlemma} we can choose a consistent pair
$(\hat{\varpi},\hat{\nu})$ such that
$d(\nu,\hat\nu)<\delta$ and, for all $a,b\in\skrix$,
\[
\sfrac{\nu_1(a)}{\hat{\nu}_1(a)}\le1+\sfrac{\delta}{\|\varpi\|+1}
\]
and
\[
\hat{\varpi}(a,b) \biggl(1+\sfrac{\delta}{\|\varpi\|+1} \biggr)\le\varpi(a,b)
\biggl(1+\sfrac{\delta}{\varpi(a,b)} \biggr).
\]
The random variables ${\ell}_j^{a}(b), j=1,\ldots,n\omega_n(a)$ are
i.i.d. with expectation
\[
\me{\ell}_1^{a}(b)=\sum_{\ell}\hat{\nu}(\ell| a)\ell(b)=\frac{ \hat
{\varpi}(a,b) } {\hat\nu_1(a)} .
\]
Hence, by the strong law of large numbers, almost surely,
\[
\limsup_{n\rightarrow\infty} \Biggl(\frac1n
\sum_{j=1}^{n\omega_n(a)}{\ell}_{j}^{a}(b)-\varpi_n(a,b) \Biggr)\le\sfrac{\nu
_1(a)}{\hat{\nu}_1(a)}\hat{\varpi}(a,b)-\varpi(a,b)\le
\delta,
\]
where we also used that $\omega_n(a)\to\omega(a)=\nu_1(a)$ and
$\varpi_n(a,b)\to\varpi(a,b)$.
\end{pf}
\begin{pf*}{Proof of Lemma \ref{randomg.approx}} We use a randomized
construction.
Given $(\varpi,{\nu})$ sub-consistent with ${\nu}_1=\omega$ and $\eps
>0$, choose
$\hat\nu$ as in Lemma \ref{random.approx1} with $\delta=\eps/(3m)$,
where $m$
is the cardinality of $\skrix$. For every $a\in\skrix$, we draw tuples
$\ell_j^{a}$,
$j=1,\ldots,n\omega_n(a)$ independently according to
$\hat\nu( \cdot| a)$ and define $e_n(a,b)$ by
\[
e_n(a,b):=\frac1n
\sum_{j=1}^{n\omega_n(a)}\ell_{j}^{a}(b)-\varpi_n(a,b)\qquad \mbox{for
all $a,b\in\skrix$}.
\]
We modify the tuples $(\ell_j^{a} \dvtx
j=1,\ldots,n\omega_n(a))$ as follows:
\begin{itemize}
\item if $e_n(a,b)<0$, we add an amount to the last element
$\ell_{ n\omega_n(a)}^{a}(b)$ such that the modified tuple satisfies
$e_n(a,b)=0$;
\item if $e_n(a,b)>0$, by Lemma \ref{random.approx1}, the ``overshoot''
$ne_n(a,b)$ cannot exceed $n\delta$. We successively deduct one
from the nonzero elements in $\ell_{ j}^{a}(b)$,
$j=1,\ldots,n\omega_n(a)$ until the modified tuples
satisfy $e_n(a,b)=0$;
\item if $e_n(a,b)=0$ we do not modify $\ell_{j}^{a}(b)$.
\end{itemize}
We denote by $(\tilde{\ell}_{j}^{a}\dvtx j=1,\ldots,n\omega_n(a))$
the tuples after all modifications.

For each $a\in\skrix$ define probability measures $\tilde{\Delta}_n(
\cdot | a)$
and $\Delta_n( \cdot| a)$ by
\[
\tilde{\Delta}_n(\ell| a)=\frac{1}{n\omega_n(a)}\sum_{j=1}^{n\omega_n(a)}
\one_{\{ \tilde{\ell}_j^{a} = \ell\}} \qquad\mbox{for
$\ell\in\skrin(\skrix)$},
\]
respectively,
\[
\Delta_n(\ell| a)=\frac{1}{n\omega_n(a)}\sum_{j=1}^{n\omega_n(a)}
\one_{\{\ell_j^{a}=\ell\}} \qquad\mbox{for $\ell\in\skrin(\skrix)$}.
\]
We define probability measures $\tilde{\nu}_n\in\skrim_n(\skrix\times
\skrin(\skrix))$
and $\nu_n\in\skrim_n(\skrix\times\skrin(\skrix))$ by
$\tilde{\nu}_n(a,\ell)=\omega_n(a){\Delta}_n(\ell| a)$, respectively,
$\nu_n(a,\ell)=\omega_n(a)\tilde{\Delta}_n(\ell| a)$, for
$(a,\ell)\in\skrix\times\skrin(\skrix)$. Recall from our modification
procedure that,
in the worst case, we have changed $n m\delta$ of the tuples. Thus,
\[
d(\tilde{\nu}_n, \nu_n) \le m\delta\le\tfrac13 \eps .
\]
As a result of our modifications
we have $\Phi(\nu_n)=(\omega_n,\varpi_n)$. We observe that, for all
$(a,\ell)\in\skrix\times\skrin(\skrix)$, the random variables
\[
\one\{\ell_1^{a}=\ell\},\ldots,\one\bigl\{\ell_{n\omega_n(a)}^{a}=\ell\bigr\}
\]
are independent Bernoulli random variables with success probability
$\hat\nu(\ell| a)$ and hence, almost surely,
\[
\lim_{n\rightarrow\infty}\Delta_n(\ell| a)=\hat\nu(\ell| a) .
\]
Therefore, for all $(a,\ell)\in\skrix\times\skrin(\skrix)$, we obtain
$\lim_{n\rightarrow\infty}\tilde{\nu}_n(a,\ell)=\hat\nu(a,\ell)$,
almost surely. Thus, almost surely, for all large $n$,
we have $d(\nu_n,\nu) \le d(\nu_n,\tilde\nu_n) + d(\tilde\nu_n, \hat\nu)
+ d(\hat\nu, \nu) \le\eps$, as claimed.
\end{pf*}
\begin{lemma}[(Approximation step 3)]\label{ntoathird}
Let $\nu_n\in\skrim_n(\skrix\times\skrin(\skrix))$ with
$\Phi(\nu_n)=(\omega_n,\varpi_n)$. For every \mbox{$\eps>0$} there
exists $n(\eps)$ such that, for all $n\ge n(\eps)$,
we can find $\tnu_n\in\skrim_n(\skrix\times\skrin(\skrix))$
satisfying \eqref{degr} and $\Phi(\tnu_n)=(\omega_n,\varpi_n)$.
\end{lemma}
\begin{pf} As $\nu_n\in\skrim_n(\skrix\times\skrin(\skrix))$, there
is a representation
\[
\nu_n=\frac1n \sum_{k=1}^n \delta_{(a_k,\ell_k)}\qquad
\mbox{for } a_k\in\skrix, \ell_k\in\skrin(\skrix).
\]
Fix $\delta>0$ and $a\in\skrix$. Look at the sets:
\begin{itemize}
\item$V^+=\{ 1\le k\le n \dvtx a_k=a, \sum_b \ell_k(b) >n^{1/3}\}$
with cardinality\vspace*{1pt} $\sharp V^+\le(n\times\break\sum_b \varpi_n(a,b))^{2/3}$,
\item$V^-=\{ 1\le k\le n \dvtx a_k=a, \sum_b \ell_k(b) \le n^{1/4}
\}$ with cardinality $\sharp V^- \ge
n- (n\sum_b \varpi_n(a,b))^{3/4}$.
\end{itemize}
For each $k\in V^+$ we replace $\ell_k$ by a smaller vector $\tilde\ell
_k$ such that
$\sum_b \tilde\ell_k(b)=n^{1/3}$. As
\[
\sum_{k\in V_1} \sum_b \ell_k(b)\le n\sum_{b\in\skrix} \varpi_n(a,b)
\]
we may replace (for large $n$) no more than $\delta n$ of the vectors
$\ell_k$, $k\in V^-$, by
larger vectors $\tilde\ell_k$ such that
\[
\sum_b \tilde\ell_k(b)\le n^{1/3}
\]
and
\[
\sum_{k=1}^n \sum_{b\in\skrix} \one\{a_k=a\} \tilde\ell_k(b)=
\sum_{k=1}^n \sum_{b\in\skrix} \one\{a_k=a\} \ell_k(b) ,
\]
where we use the convention $\tilde\ell_k=\ell_k$ if this vector was
not changed in
the procedure. Performing such an operation for every $a\in\skrix$ we
may define
\[
\tilde\nu_n=\frac1n \sum_{k=1}^n \delta_{(a_k,\tilde\ell_k)},
\]
and observe that \eqref{degr} holds and $\Phi(\tnu_n)=(\omega_n,\varpi_n)$.
Moreover,
\[
d(\nu_n,\tnu_n) \le\sfrac{m}{2n} \biggl( \biggl(n\sum_b \varpi_n(a,b)\biggr)^{2/3} +
\delta n \biggr),
\]
which is less than $\eps>0$ for a suitable choice of $\delta>0$, and
all sufficiently large $n$.
\end{pf}

\subsection{\texorpdfstring{Proof of the lower bound in Theorem \protect\ref
{main3}}{Proof of the lower bound in Theorem 2.5}}\label{sec45}

There is a partial analogue to Lemma \ref{upperb} for the lower bounds.
\begin{lemma}\label{lowerb}
For any $\nu_n\in\skrim_n(\skrix\times\skrin(\skrix))$ which satisfies
\eqref{degr}
with $\Phi(\nu_n)=(\omega_n,\varpi_n)$ and any $\eps>0$, we have
\[
{\prob} \{ d(M,\nu_n)<\eps | \Phi(M)=(\omega_n,\varpi_n) \}
\ge\exp\bigl( -nH(\nu_n \| Q_n)-\eps_{2}^{\ssup n}(\nu_n) \bigr) ,
\]
where $Q_n$ is as Lemma \ref{upperb} and
\[
\lim_{n\uparrow\infty} \frac1n \eps_{2}^{\ssup n}(\nu_n)=0 .
\]
\end{lemma}
\begin{pf*}{Proof of Lemma \ref{lowerb}}
We use the notation and some results from the proof of the upper bound,
Lemma \ref{upperb}. In particular, recall
the definition of $n(a,b)$ from \eqref{nabdef} and, from \eqref
{randomg.probcom}, that
%
%
\begin{eqnarray}\label{lbuni}
&&{\prob} \{d(M,\nu_n)<\eps | \Phi({M})=(\omega_n,\varpi_n) \}\nonumber\\[-8pt]\\[-8pt]
&&\qquad=\frac{\sharp\{\tilde{x}\in\skria_n(\omega
_n,\varpi_n) \dvtx
d({M \circ\Psi}(\tilde{x}),\nu_n)< \eps\}}{\sharp\{\tilde{x}\in
\skria_n(\omega_n,\varpi_n) \}},\nonumber
\end{eqnarray}
and that the denominator was evaluated in \eqref{denom} as
\[
\pmatrix{n\cr n\omega_n(a), a\in\skrix}
\prod_{\{a,b\}} \prod_{k=1}^{n(a,b)}
\biggl( \frac{n^2\omega_n(a)\omega_n(b) - n\omega_n(a)\one_{\{a=b\}}}{1+\one
_{\{a=b\}}} - (k-1) \biggr) .
\]
We now describe a procedure which yields (for sufficiently large $n$) a
lower bound of
%
%
\begin{eqnarray}\label{numlb}
&&\pmatrix{n\cr n\nu_n(a,\ell), a\in\skrix, \ell\in\skrin(\skrix)}
\nonumber\\[-8pt]\\[-8pt]
&&\qquad{}\times\prod_{(a,b)} \frac{ (n\varpi_n(a,b)-2\lceil n^{2/3} \rceil-2)!}{ \prod
_{j=1}^{n\omega_n(a)}
(\ell_{a}^{\ssup j}(b))!} 2^{-n/2 \varpi_n(\Delta)} \nonumber
\end{eqnarray}
for the numerator, where $\ell_{a}^{\ssup j}(b), j=1,\ldots,n\omega
_n(a)$ are any enumeration of the family
containing each $\ell(b)$ with multiplicity $n\nu_n(a,\ell)$.

First, we allocate to each vertex some
$(a,\ell)\in\skrix\times\skrin(\skrix)$ in such a way that every
vector $(a,\ell)$ is allocated $n\nu_n(a,\ell)$ times.
There are
\[
\pmatrix{n\cr n\nu_n(a,\ell), a\in\skrix, \ell\in\skrin(\skrix)}
\]
ways to do this.
Next we add edges between the vertices of two \textit{different} colors,
$a, b\in\skrix$.
To this end, we distribute the numbers in $\{1,\ldots,n(a,b)-\lceil
n^{2/3}\rceil-1\}$
among the vertices with color $a$ so that a vertex carrying vector
$(a,\ell)$ gets at most $\ell(b)$ numbers.
A crude lower bound for the number of ways to do this is
\[
\frac{ (n(a,b)-\lceil n^{2/3} \rceil-1 )!}{ \prod_{j=1}^{n\omega_n(a)}
(\ell_{a}^{\ssup j}(b))!}.
\]
Now the numbers in $\{1,\ldots,n(a,b)-\lceil n^{2/3}\rceil-1\}$ are
distributed successively, this time among the vertices
of color $b$. Again we do this in such a way that a vertex carrying
vector $(b,\ell)$ has a capacity to carry no
more than $\ell(a)$ numbers. However, we are more cautious now: when
distributing $k$ we look at the vertex of color $a$,
which already carries $k$. If this carries numbers from $\{1,\ldots
,k-1\}$, we do not allow $k$ to be associated with any
vertex of color $b$ which carries one of these numbers.
By \eqref{degr} this rules out no more than $n^{1/3}$ vertices, each of which
has a capacity no more than $n^{1/3}$, so that the number of ways to do
this is at least
\[
\frac{ (n(a,b)-\lceil n^{2/3} \rceil-1)!}{ \prod_{j=1}^{n\omega
_n(b)}(\ell_{b}^{\ssup j}(a))!}.
\]
Next, for $k=1,\ldots, n(a,b)-\lceil n^{2/3} \rceil-1$ draw the $k$th
edge between the two vertices of color $a$ and $b$
carrying number $k$ and observe that when allocating the numbers we
have been cautious not to cause any multiple edges.
Obviously, there is at least one way to establish a further $\lceil
n^{2/3} \rceil+1$ edges between vertices of
color $a$ and $b$ without creating multiple edges.

We now add the edges connecting vertices of \textit{the same} color
$a\in
\skrix$. For this purpose, we successively
distribute the numbers in $\{1,\ldots,n(a,a)-\lceil n^{2/3}\rceil-1\}$
and $\{n(a,a)+1,\ldots,2n(a,a)-\lceil n^{2/3}\rceil-1\}$
among the vertices with color $a$ so that a vertex carrying vector
$(a,\ell)$ gets at most $\ell(b)$ numbers.
When distributing $k>n(a,a)$ we look at the vertex of color $a$, which
already carries $k-n(a,a)$.
If this vertex carries numbers $j\in\{1,\ldots,k-n(a,a)-1\}$, we do not
allow $k$ to be associated with the vertices
carrying numbers $j+n(a,a)$. We also do not allow $k$ to be associated
with the vertex itself.
By \eqref{degr} these restrictions rule out no more than $n^{1/3}$
vertices, each of which
has a capacity no more than $n^{1/3}$, so that the number of ways to do
this is at least
\[
\frac{ (2n(a,a)-2\lceil n^{2/3} \rceil-2)!}{ \prod_{j=1}^{n\omega
_n(a)}(\ell_{a}^{\ssup j}(a))!}.
\]
Obviously, there is at least one way of allocating the remaining
numbers $\{n(a,a)-\lceil n^{2/3}\rceil, \ldots, n(a,a)\}$ and
$\{2n(a,a)-\lceil n^{2/3}\rceil, \ldots, 2n(a,a)\}$ to vertices so that
no single vertex carries a matching pair $j, j+n(a,a)$,
and no pair of vertices carry two or more matching pairs between them.
Next, for $k=1,\ldots, n(a,b)$ draw the $k$th edge between
the two vertices carrying numbers $k$ and $k+n(a,a)$ and observe that
when allocating the numbers we have been cautious not to cause
any loops or multiple edges. As, for every $k\in\{1,\ldots, n(a,a)\}$,
the numbers $k$ and $k+n(a,a)$ could be interchanged
without changing the configuration, the total number of different
configurations constructable in this procedure is bounded from below by
\begin{eqnarray*}
&&\pmatrix{n\cr\nu_n(a,\ell), a\in\skrix, \ell\in\skrin(\skrix)}
\times
\mathop{\prod_{(a,b)}}_{a\not= b} \frac{ (n(a,b)-\lceil n^{2/3} \rceil
-1 )!}{ \prod_{j=1}^{n\omega_n(a)}
(\ell_{a}^{\ssup j}(b))!} \\
&&\qquad{}\times\prod_{a\in\skrix} \frac{
(2n(a,a)-2\lceil n^{2/3} \rceil-2 )!}{ \prod_{j=1}^{n\omega_n(a)}
(\ell_{a}^{\ssup j}(a))!} 2^{-n(a,a)} ,
\end{eqnarray*}
and this is bounded from below by the quantity in \eqref{numlb}. Every
resulting graph satisfies the constraint
$\Phi(M)=(\omega_n, \varpi_n)$. To measure the distance between its
empirical neighborhood measure $M$
and $\nu_n$, we say that a vertex $v\in V$ is \textit{successful} if the
associated $(X(v),L(v))$ is identical
to the $(a,\ell)$ they were carrying after the initial step. Note that
after\vspace*{1pt} allocation of the edges with numbers
in $\{1,\ldots,n(a,b)-\lceil n^{2/3}\rceil-1\}$ among the vertices of
all colors, all but at most
$2m^2(\lceil n^{2/3}\rceil+1)$ vertices $v\in V$ were successful.
Adding in the further edges in the last step
can lead to up to $2m^2(\lceil n^{2/3}\rceil+1)$ further unsuccessful
vertices. Hence
\[
d(\nu_n,M) \le\frac1{2n} \sum_{v\in V} \one\{v \mbox{ unsuccessful}\}
\le4m^2 n^{-1/3}
\stackrel{n\to\infty}{\longrightarrow} 0\qquad
\mbox{as }n\uparrow\infty.
\]
To complete the proof, we again use Stirling's formula to analyze the
combinatorial terms obtained as an estimate
for the numerator and denominator in \eqref{lbuni}. 
For the denominator we get
the same main terms as in Lemma \ref{upperb} with slightly different
error terms, which however do not
depend on $\nu_n$. More interestingly, we have
\begin{eqnarray*}
&&
\prod_{a\in\skrix} \pmatrix{n\omega_n(a)\cr n\nu_n(a,\ell), \ell\in
\skrin(\skrix)}\\[-1pt]
&&\qquad \ge\exp\biggl(n\sum_{a}\omega_n(a)\log\omega_n(a)-n\sum_{(a,\ell)}\nu
_n(a,\ell)
\log\nu_n(a,\ell) \biggr)\\[-1pt]
&&\qquad\quad{} \times\exp\biggl( -\sfrac{|\skris(\nu_n)|}{2}\log(2\pi n)-
\mathop{\sum_{(a,\ell)}}_{n\nu_n(a,\ell)\ge1}\sfrac{1}{12n\nu_n(a, \ell
)} \biggr) ,
\end{eqnarray*}
where the exponent in the error term is of order o$(n)$, by the bound
on the size of the support of $\nu_n$
given in Lemma \ref{randomg.finiteness}. Further,
\begin{eqnarray*}
\hspace*{-4pt}&&\frac{ (n\varpi_n(a,b)-2\lceil n^{2/3} \rceil-2)!}{ \prod
_{j=1}^{n\omega_n(a)}
(\ell_{a}^{\ssup j}(b))!}\\
\hspace*{-4pt}&&\qquad \ge\exp\biggl( -n\sum_{\ell} \log(\ell(b)! ) \nu_n(a,\ell)
+n\varpi_n(a,b)\log(n\varpi_n(a,b) )-n\varpi_n(a,b) \biggr) \\
\hspace*{-4pt}&&\qquad\quad{} \times\exp\biggl(-(2n^{2/3}+2) \log(n\varpi_n(a,b) )
+ n\varpi_n(a,b) \log\biggl(1- \sfrac{2n^{2/3}+2}{n\varpi_n(a,b)} \biggr)\biggr) ,
\end{eqnarray*}
and the result follows by combining this with facts discussed in the
context of the upper
bound.
\end{pf*}

To complete the proof of the lower bound in Theorem \ref{main3},
take an open set $\Gamma\subset\skrim(\skrix\times\skrin(\skrix))$.
Then, for any $\nu\in\Gamma$ with $(\varpi,{\nu})$ sub-consistent and
${\nu}_1=\omega$ we may find $\eps>0$ with the ball around $\nu$ of
radius $2\eps>0$ contained in $\Gamma$. By our approximation (see
Lemmas \ref{randomg.approx}
and \ref{ntoathird}) we may find $\nu_n\in\Gamma\cap\skrim_n(\skrix
\times\skrin(\skrix))$ with
$\Phi(\nu_n)=(\omega_n,\varpi_n)$ such that \eqref{degr} holds and
$d(\nu_n,\nu)\downarrow0$.
Hence, for all large $n\ge n(\eps)$,
\begin{eqnarray*}
{\prob} \{M\in\Gamma | \Phi(M)=(\omega_n,\varpi_n) \}
& \ge & {\prob} \{ d(\nu_n, M)< \eps | \Phi(M)=(\omega_n,\varpi_n) \}
\\
& \ge &\exp\bigl( -nH(\nu_n \| Q_n)-\eps_2^{\ssup n}(\nu_n) \bigr).
\end{eqnarray*}
We observe that
\begin{eqnarray*}
\lim_{n\rightarrow\infty} H(\nu_n \| Q_n) - H(\nu\| Q) &=& \lim
_{n\rightarrow\infty}
H(\nu_n \| Q_n) - H(\nu_n \| Q)\\
&&{} +\lim_{n\rightarrow\infty}H(\nu_n \|
Q) - H(\nu\| Q)\\
&=&0,
\end{eqnarray*}
where the last term vanishes by continuity of relative entropy, and the
first term was shown to vanish
in the proof of Lemma \ref{lowerb}. This completes the proof of
Theorem \ref{randomg.LDprob}.

\section{\texorpdfstring{Proof of Theorems \protect\ref{main2} and \protect\ref
{main}(b) by mixing}{Proof of Theorems 2.1 and 2.3(b) by mixing}}\label{proofmain2}

We denote by $\Theta_n:=\skrim_n(\skrix)\times\tilde\skrim_{*,n}(\skrix
\times\skrix)$
and $\Theta:=\skrim(\skrix)\times\tilde{\skrim}_*(\skrix\times\skrix)$. Define
\begin{eqnarray*}
P_{(\omega_n, \varpi_n)}^{\ssup n}(\nu_n) & := &\prob\{M=\nu_n |
\Phi(M)=(\omega_n,\varpi_n) \} ,\\
P^{\ssup n}(\omega_n,\varpi_n) & := &
\prob\{(L^1,L^2)=(\omega_n,\varpi_n) \} ,\\
P_{\omega_n}^{\ssup n}(\varpi_n)
& := &\P\{L^2=\varpi_n | L^1=\omega_n\} ,\\
P^{\ssup n}(\omega_n) & := & \P\{L^1=\omega_n\} .
\end{eqnarray*}
The joint distribution of $L^1, L^2$ and $M$ is the mixture of
$P_{(\omega_n, \varpi_n)}^{\ssup n}$ with $P^{\ssup n}(\omega_n$, $\varpi
_n)$ defined as
%
%
\begin{equation}\label{randomg.mixture}
d\tilde{P}^n(\omega_n, \varpi_n, \nu_n):= dP_{(\omega_n,
\varpi_n)}^{\ssup n}(\nu_n) \,dP^{\ssup n}(\omega_n, \varpi_n) ,
\end{equation}
while the joint distribution of $L^1$ and $L^2$ is the mixture of
$P_{\omega_n}^{\ssup n}$
with $ P^{\ssup n}$ given by
%
%
\begin{equation}\label{randomg.mixture1}
dP^{\ssup n}(\omega_n, \varpi_n)=dP_{\omega_n}^{\ssup n}(\varpi_n)\,
dP^{\ssup n}(\omega_n).
\end{equation}

Biggins (\cite{Bi04}, Theorem 5(b)) gives criteria for the
validity of large deviation principles for the mixtures and for the
goodness of the rate function if individual large deviation
principles are known. The following two lemmas ensure validity of
these conditions.
\begin{lemma}[(Exponential tightness)]\label{randomg.tightness}
The following families of distributions are exponentially tight:
\begin{enumerate}[(a)]
\item[(a)] $(P^{\ssup n} \dvtx n\in\N)$ on $\skrim(\skrix)\times
\tilde{\skrim}_*(\skrix\times\skrix)$,
\item[(b)] $(\tilde{P}^{\ssup n} \dvtx n\in\N)$ on $\skrim(\skrix
)\times\tilde{\skrim}_*(\skrix\times\skrix)\times\skrim(\skrix\times
\skrin(\skrix))$.
\end{enumerate}
\end{lemma}
\begin{pf}
(a) It suffices to show that, for every $\theta>0$, there exists $N\in
\N$ such that
\[
\limsup_{n\rightarrow\infty}\sfrac{1}{n}\log\prob\{|E|>n N \}\le-\theta.
\]

To see this, let $c>\max_{a,b\in\skrix} C(a,b)>0$. By a simple coupling argument
we can define, for all sufficiently large $n$, a new colored random
graph $\tilde{X}$ with color law $\mu$ and connection probability
$\frac{c}{n}$, such that any edge present in $X$ is also present in
$\tilde{X}$. Let $|\tilde{E}|$ be the number of edges of
$\tilde{X}$. Using Chebyshev's inequality, the binomial formula and
Euler's formula, we have that
\begin{eqnarray*}
\prob\{ |\tilde{E}|\ge n l \}
&\le&
e^{-nl}\me\bigl[e^{|\tilde{E}|} \bigr]\\
&=& e^{-nl}\sum_{k=0}^{{n(n-1)}/{2}}e^{k} \pmatrix{n(n-1)/2\cr
k} \biggl(\frac{c}{n} \biggr)^{k} \biggl(1-\frac{c}{n} \biggr)^{n(n-1)/2-k}\\
&=& e^{-nl} \biggl( 1-\frac{c}{n}+e\frac{c}{n} \biggr)^{n(n-1)/2} \\
&\le&
e^{-nl}e^{nc(e-1+o(1))}.
\end{eqnarray*}
Now given $\theta>0$ choose $N\in\N$ such that $N> \theta+c(e-1)$
and observe that, for sufficiently large $n$,
\[
\prob\{|E|\ge n N \}\le\prob\{ |\tilde{E}|\ge
nN \}\le e^{-n\theta},
\]
which implies the statement.

(b) Given $\theta>0$, we observe from (a) that there exists $N(\theta
)\in\N$ such that,
for all sufficiently large $n$,
\[
\prob\bigl\{M\bigl(\{\|\ell\|\ge2\theta N(\theta)\}\bigr)\ge\theta^{-1} \mbox{ or }
\|L^2\| \ge2N(\theta) \bigr\}
\le\prob\{|E|\ge nN(\theta) \}\le e^{-\theta n}.
\]
We define the set
$\Xi_{\theta}$ by
\begin{eqnarray*}
&&\Xi_{\theta} : = \bigl\{ (\varpi,\nu)\in\tilde\skrim_*(\skrix\times\skrix
)\times\skrim\bigl(\skrix\times\skrin(\skrix)\bigr)
\dvtx\nu\{\|\ell\|>2lN(l)\}< l^{-1}\\
&&\hspace*{153.5pt}\hspace*{56.9pt} \forall l\ge\theta\mbox{ and }
\|\varpi\|<2N(\theta) \bigr\}.
\end{eqnarray*}
As $\{\|\ell\|\le2l N(l)\}\subset\skrin(\skrix)$ is finite, hence compact,
the set $\Xi_\theta$ is relatively compact in the weak topology, by
Prohorov's criterion. Moreover, we have that
\[
\tilde{P}^{n} ((\Xi_\theta)^{\mathrm{c}} )\le\prob\{ \|L^2\| \ge
2N(\theta
) \} +
\sum_{l=\theta}^{\infty} \prob\bigl\{ M\bigl(\{\|\ell\|>2lN(l)\}\bigr)\ge
l^{-1} \bigr\} \le C(\theta) e^{-n\theta}.
\]
Therefore, $ \limsup_{n\rightarrow\infty}\sfrac1n\log
\tilde{P}^{n}((\operatorname{cl} \Xi_\theta)^{\mathrm{c}})\le-\theta$, which
completes the proof, as $\theta>0$ was arbitrary.
\end{pf}

Now, we observe that the function
$I(\omega,\varpi)=H(\omega\| \mu)+ \skrih_C(\varpi\|\omega)$ is a
good rate function,
by a similar argument as in the proof of Lemma \ref{randomg.Vrate}(ii).
Therefore, applying Theorem 5(b) of \cite{Bi04} to Sanov's
theorem (\cite{DZ98}, Theorems 2.1.10 and~\ref{main}(a)), we obtain the large deviation principle for
$P^{\ssup n}$ on $\skrim(\skrix)\times\tilde{\skrim}_{*}(\skrix\times
\skrix)$ with
good rate function $I(\omega,\varpi)$, which is Theorem \ref{main}(b).

To prove Theorem \ref{main3} define the function
\[
\tilde{J}\dvtx{\Theta}\times\skrim\bigl(\skrix\times\skrin(\skrix
)\bigr)\rightarrow[0,\infty],\qquad
\tilde{J}((\omega,\varpi), \nu)=\tilde{J}_{(\omega,\varpi)}(\nu).
\]

\begin{lemma}\label{randomg.convexgoodrate}
$\tilde{J}$ is lower semicontinuous.
\end{lemma}
\begin{pf} Suppose $\theta_n:=((\omega_n,\varpi_n), \nu_n)$ converges to
$\theta:=((\omega,\varpi), \nu)$ in ${\Theta}\times\skrim(\skrix\times
\skrin(\skrix))$.
There is nothing to show if $\liminf_{\theta_{n}\rightarrow\theta}
\tilde J(\theta_n)=\infty$.
Otherwise, if $(\varpi_n, \nu_n)$ is sub-consistent for infinitely many
$n$, then
\[
\varpi(a,b)= \lim_{n\uparrow\infty} \varpi_n(a,b) \ge
\liminf_{n\uparrow\infty} \langle\nu_n(\cdot,\ell), \ell(\cdot)
\rangle(a,b)
\ge\langle\nu(\cdot,\ell), \ell(\cdot) \rangle(a,b) ,
\]
hence $(\varpi, \nu)$ is sub-consistent. Similarly, if the first marginal
of $\nu_n$ is $\omega_n$, we see that the first marginal of $\nu$ is
$\omega$.
We may therefore argue as in \eqref{equ.entropyproof} to obtain
\begin{eqnarray*}
\liminf_{\theta_{n}\rightarrow\theta} J(\theta_n) &=&
\liminf_{\theta_{n}\rightarrow\theta} H(\nu_n \| Q_n)\\
&\ge&
\lim_{\theta_{n}\rightarrow\theta} H(\nu_n \| Q_n)
-H(\nu_n \| Q) + \liminf_{\nu_{n}\rightarrow\nu}H(\nu_n \| Q)\\
&=&H(\nu\| Q),
\end{eqnarray*}
where the last step uses continuity
of relative entropy. This proves the lemma.
\end{pf}

Lemmas \ref{randomg.convexgoodrate} and \ref{randomg.tightness}(b)
ensure that we can apply Theorem 5(b) of \cite{Bi04} to the large deviation
principles established in
Theorem \ref{main}(b) and \ref{main3}. This yields a large deviation
principle for
$(\tilde{P}^n \dvtx n\in\N)$ on $\skrim(\skrix) \times\tilde\skrim
_*(\skrix\times\skrix)
\times\skrim(\skrix\times\skrin(\skrix))$ with good rate function
\[
\hat{J}(\omega, \varpi, \nu)= \cases{H(\omega\| \mu) + \sfrac12
{\mathfrak H}_C(\varpi\| \omega) +
H(\nu\| Q) , &\quad if $(\varpi,{\nu})$ sub-consistent\cr
&\quad and ${\nu}_1=\omega$,\cr
\infty , &\quad otherwise.}
\]
By projection onto the last two components we obtain
the large deviation principle as stated in Theorem \ref{main2}
from the contraction principle; see, for example, \cite{DZ98}, Theorem 4.2.1.

\section{\texorpdfstring{Proof of Corollaries \protect\ref{ERdd} and
\protect\ref{randomg.L1E}}{Proof of Corollaries 2.2 and 2.4}}
\label{proofcorollaries}

We derive the corollaries from Theorem~\ref{main2} by applying the
contraction principle;
see, for example, \cite{DZ98}, Theorem~4.2.1. It just remains to
simplify the rate functions.

\subsection{\texorpdfstring{Proof of Corollary \protect\ref{ERdd}}{Proof of Corollary 2.2}}

In the case of an uncolored Erd\H{o}s--Renyi graph, the function $C$
degenerates to a constant $c$,
$L^2=|E|/n\in[0,\infty)$ and $M=D\in\skrim(\N\cup\{0\})$. Theorem \ref
{main2} and the
contraction principle imply a large deviation principle for $D$ with
good rate function
\begin{eqnarray*}
\delta(d)&=&\inf\{J(x,d) \dvtx x\ge0 \}\\
&=&\inf\bigl\{ H(d \| q_{x})+\tfrac1 2 x\log
x -\tfrac12 x\log c+ \tfrac{1}{2} c-\tfrac12 x \dvtx\langle d\rangle
\le x \bigr\},
\end{eqnarray*}
which is to be understood as infinity if $\langle d\rangle$ is
infinite. We denote by $\delta^{x}(d)$ the expression inside the
infimum and consider the cases (i) $\langle d\rangle\le c$ and
(ii) $\infty>\langle d\rangle\ge c$ separately.

\textit{Case} (i). Under our condition, the equation
$x=ce^{-2(1-\langle d\rangle/x)}$
has a unique solution, which satisfies $x\ge\langle d\rangle$.
Elementary calculus shows that the global minimum of
$y\mapsto\delta^{y}(d)$ on $(0,\infty)$ is attained at the value $y=x$,
where $x$ is the solution of our equation.

\textit{Case} (ii). For any $\eps>0$, we have
\begin{eqnarray*}
\delta^{\langle d\rangle+\eps}(d)-\delta^{\langle d\rangle}(d)
&=& \sfrac{\eps}{2} +\sfrac{\langle d\rangle-\eps}{2}\log\sfrac{\langle
d\rangle}{\langle
d\rangle+\eps}+\sfrac{\eps}{2}\log\sfrac{\langle d\rangle}{c}\\
&\ge&\sfrac{\eps}{2} +\sfrac{\langle d\rangle-\eps}{2} \biggl(\sfrac{-\eps
}{\langle
d\rangle} \biggr)+\sfrac{\eps}{2}\log\sfrac{\langle d\rangle}{c} >0,
\end{eqnarray*}
so that the minimum is attained at $x=\langle d\rangle$.

\subsection{\texorpdfstring{Proof of Corollary \protect\ref{randomg.L1E}}{Proof of Corollary 2.4}}
We begin the proof by defining the continuous linear map
$W\dvtx\skrim(\skrix)\times\tilde{\skrim}_*(\skrix\times\skrix)
\rightarrow[0,\infty)$ by $W(\omega,{\varpi})=\sfrac{1}{2}\|{\varpi}\|$.
We infer from Theorem \ref{main} and the contraction
principle that $W(L^1,L^2)=|E|/n$
satisfies a large deviation principle in $[0,\infty)$ with the good
rate function
\[
\zeta(x)=\inf\{I(\omega,{\varpi}) \dvtx W(\omega,{\varpi})=x \}.
\]
To obtain the form of the rate in the corollary, the infimum is
reformulated as unconstrained
optimization problem (by normalizing $\varpi$)
%
%
\begin{equation}\label{randomg.Jen}
\mathop{\inf_{\varpi\in\skrim_*(\skrix\times\skrix)}}_{\omega\in\skrim
(\skrix)}
\biggl\{H(\omega\| \mu)+xH(\varpi\| C\omega\otimes\omega)+x\log2x+
\sfrac{1}{2}
\| C \omega\otimes\omega\| -x \biggr\}.\hspace*{-28pt}
\end{equation}
By Jensen's\vspace*{1pt} inequality $H(\varpi\| C\omega\otimes\omega)\ge-\log\|
C\omega\otimes\omega\|$,
with equality if $\varpi=\sfrac{C\omega\otimes\omega}{\| C\omega\otimes
\omega\|}$, and hence,
by symmetry of $C$ we have
\begin{eqnarray*}
&&\min_{\varpi\in\skrim_*(\skrix\times\skrix)} \biggl\{ H(\omega\| \mu
)+xH(\varpi\| C\omega\otimes\omega)
+x\log2x+\sfrac12 \| C \omega\otimes\omega\|-x \biggr\}\\
&&\qquad= H(\omega\| \mu) -x\log\| C \omega\otimes\omega\|
+x\log2x + \sfrac12 \| C \omega\otimes\omega\|-x.
\end{eqnarray*}

The form given in Corollary \ref{randomg.L1E} follows by defining
\[
y=\sum_{a,b\in\skrix} C(a,b)\omega(a)\omega(b).
\]

\section*{Acknowledgments}
We thank Nadia Sidorova for helpful discussions and, in particular, for
contributing the idea of the proof of Lemma \ref{randomg.dpartition};
Peter Eichelsbacher for helpful discussions on large deviations
principles in strong topolgies; and two referees for numerous
suggestions which significantly improved the paper. In particular, the
proof of Theorem \ref{main} was considerably shortened by a referee's
suggestion to use G\"artner--Ellis and mixing instead of a direct
change-of-measure argument. This paper contains material from the first
author's Ph.D. thesis.

Revision of this paper was done during a visit of the first author to
Bath, which was supported by the London Mathematical Society.

%

%
\printaddresses

\end{document}